\documentclass{amsart}
\usepackage[dvips]{graphicx}
\usepackage{ifthen,amsfonts,amsmath,amssymb,epic,eepic}  %put in
\usepackage{epsfig,color}  %put in

\theoremstyle{plain}
\newtheorem{theorem}{Theorem}[section]
\newtheorem{lemma}[theorem]{Lemma}
\newtheorem{corollary}[theorem]{Corollary}
\newtheorem{proposition}[theorem]{Proposition}

\newtheorem{conjecture}{Conjecture}

\theoremstyle{remark}

\newtheorem{remark}[theorem]{\bf Remark}
\theoremstyle{definition}

\newcommand{\dist}{{\text{dist}}}

\newcommand{\Energy}{{\text{Energy}}}

\newcommand{\Graph}{\operatorname{Graph}}

\def\ZZ{{\bold Z}}
\def\RP{{\bold RP}}
\def\RR{{\bold R}}
\def\SS{{\bold S}}

\newcommand{\nn}{{\bold {n}}}
\newcommand{\HH}{{\bold {H}}}

\newcommand{\dv}{{\text{div }}}

\newcommand{\Vol}{{\text {Vol}}}
\newcommand{\Area}{{\text {Area}}}
\newcommand{\Ric}{{\text {Ric}}}
\newcommand{\Scal}{{\text {Scal}}}

\def\CC{{\bold C }}
\newcommand{\eqr}[1]{(\ref{#1})}

\newcommand{\e}{{\text {e}}}
\newcommand{\cB}{{\mathcal {B}}}
\newcommand{\cP}{{\mathcal {P}}}
\newcommand{\cS}{{\mathcal {S}}}
\newcommand{\cSt}{{\mathcal{S}_{neck}}}
\newcommand{\cSu}{{\mathcal{S}_{ulsc}}}

\newcommand{\cI}{{\mathcal {I}}}
\newcommand{\cL}{{\mathcal {L}}}

\newcommand{\cM}{{\mathcal {M}}}
\newcommand{\cF}{{\mathcal {F}}}
\newcommand{\cT}{{\mathcal {T}}}
\newcommand{\Conv}{\operatorname{Conv}}

\numberwithin{equation}{section}

\begin{document}

\title[Minimal submanifolds]
{Minimal submanifolds}

\author{Tobias H. Colding}%
\address{Department of Mathematics\\
MIT\\ 77 Mass. Ave.\\ Cambridge, MA 02139}
\author{William P. Minicozzi II}%
\address{Department of Mathematics\\
Johns Hopkins University\\
3400 N. Charles St.\\
Baltimore, MD 21218}
\thanks{The authors were partially supported by NSF Grants DMS
0104453 and DMS  0405695}

\email{colding@cims.nyu.edu and minicozz@math.jhu.edu}

%\subjclass{Primary 54C40, 14E20; Secondary 46E25, 20C20}

%\date{January 1, 1994 and, in revised form, June 22, 1994.}


%\keywords{Differential geometry, algebraic geometry}

%\begin{abstract}
%\end{abstract}

\maketitle \tableofcontents

\section{Introduction}

In this article, we survey some old and new results about minimal
surfaces and submanifolds.  The field of  minimal surfaces has its
origin in the mid eighteenth century with the work of Euler and
Lagrange but it has very recently seen major advances that have
solved many long standing open conjectures in the field.  In what
follows, we give a quick tour through the field, starting with the
definition and the classical results and ending up with current
areas of research.   Many references are given for further
reading.

Soap films, soap bubbles, and surface tension were extensively
studied by the Belgian physicist and inventor (the inventor of the
stroboscope) Joseph Plateau in the first half of the nineteenth
century.  At least since his studies, it has been known that the
right mathematical model for soap films are minimal surfaces --
the soap film is in a state of minimum energy when it is covering
the least possible amount of area.   Minimal surfaces and
equations like the minimal surface equation have served as
mathematical models for many physical problems.

The field of minimal surfaces dates back to the publication in
1762 of Lagrange's famous memoir ``Essai d'une nouvelle m\'ethode
pour d\'eterminer les maxima et les minima des formules
int\'egrales ind\'efinies''.  Euler had already in a paper
published in 1744 discussed minimizing properties of the surface
now known as the catenoid, but he only considered variations
within a certain class of surfaces. In the almost one quarter of a
millennium that has past since Lagrange's memoir minimal surfaces
has remained a vibrant area of research and there are many reasons
why.  The study of minimal surfaces was the birthplace of
regularity theory.  It lies on the intersection of nonlinear
elliptic PDE, geometry,   topology and general relativity.

\part{Classical and almost classical results}       \label{p:1}

Let $\Sigma \subset \RR^n$ be a \underline{smooth} $k$-dimensional
submanifold (possibly with boundary) and $C^{\infty}_0(N\Sigma)$
 the space of all infinitely differentiable,  compactly
supported, normal vector fields on $\Sigma$.  Given $\Phi$ in
$C^{\infty}_0(N\Sigma)$, consider the one--parameter variation
\begin{equation}
\Sigma_{t,\Phi}=\{x+t\,\Phi (x)  | x\in \Sigma\}\, .
\end{equation}
The so called {\it first variation formula} of volume is the
equation
\begin{equation}  \label{e:frstvar}
\left.\frac{d}{dt} \right|_{t=0}\Vol (\Sigma_{t,\Phi})
=\int_{\Sigma} \langle \Phi\, , \, \HH \rangle \, ,
\end{equation}
where $\HH$ is the mean curvature (vector) of $\Sigma$. Here, and
throughout this paper, integration is with respect to
$d\text{vol}$.  (When $\Sigma$ is noncompact, then
$\Sigma_{t,\Phi}$ in \eqr{e:frstvar} is replaced by
$\Gamma_{t,\Phi}$, where $\Gamma$ is any compact set containing
the support of $\Phi$.) The submanifold $\Sigma$ is said to be a
{\it minimal submanifold} (or just minimal) if
\begin{equation}
\left.\frac{d}{dt} \right|_{t=0}\Vol (\Sigma_{t,\Phi})=0
\,\,\,\,\,\,\,\,\,\,\,\text{ for all } \Phi\in
C^{\infty}_0(N\Sigma)
\end{equation}
 or, equivalently by \eqr{e:frstvar}, if the
mean curvature $\HH$ is identically zero.  Thus $\Sigma$ is
minimal if and only if it is a critical point for the volume
functional. (Since a critical point is not necessarily a minimum
the term ``minimal'' is misleading, but it is time honored.  The
equation for a critical point is also sometimes called the
Euler--Lagrange equation.)

Suppose now, for simplicity, that $\Sigma$ is an oriented
hypersurface with unit normal $\nn_{\Sigma}$.  We can then write a
normal vector field $\Phi \in C^{\infty}_0(N\Sigma)$  as $\Phi =
\phi \nn_{\Sigma}$, where the function $\phi$ is in the space
$C^{\infty}_0(\Sigma)$ of infinitely differentiable, compactly
supported functions on $\Sigma$. Using this, a computation shows
that if $\Sigma$ is minimal, then
\begin{equation}    \label{e:secvar}
\left. \frac{d^2}{dt^2} \right|_{t=0}\Vol (\Sigma_{t,\phi
\nn_{\Sigma} }) =-\int_{\Sigma}\phi\,L_{\Sigma}\phi\, ,
\end{equation}
where
\begin{equation}
    L_{\Sigma}\phi=\Delta_{\Sigma}\phi+|A|^2\phi
\end{equation}
is the second variational (or Jacobi) operator. Here
$\Delta_{\Sigma}$ is the Laplacian on $\Sigma$ and $A$ is the
second fundamental form.  So $|A|^2=\kappa_1^2+\kappa_2^2+ \dots +
\kappa_{n-1}^2 $,
 where $\kappa_1,\dots \kappa_{n-1}$ are the principal curvatures of $\Sigma$
and $\HH= (\kappa_1+ \dots + \kappa_{n-1}) \, \nn_{\Sigma} $. A
minimal submanifold $\Sigma$ is said to be {\it stable} if
\begin{equation}
\left. \frac{d^2}{dt^2} \right|_{t=0}\Vol (\Sigma_{t,\Phi})\geq 0
\,\,\,\,\,\,\,\,\,\,\,\text{ for all } \Phi\in
C^{\infty}_0(N\Sigma)\, .
\end{equation}
Integrating by parts in \eqr{e:secvar}, we see that stability is
equivalent to the so called {\it stability inequality}
\begin{equation}   \label{e:stabin}
      \int |A|^2 \, \phi^2 \leq   \int |\nabla \phi|^2 \, .
\end{equation}
More generally, the {\it Morse index} of a minimal submanifold is
defined to be the number of negative eigenvalues of the operator
$L$.  Thus, a stable submanifold has Morse index zero.

\subsection{The Gauss map}

Let $\Sigma^2\subset \RR^3$ be a surface (not necessarily
minimal). The {\it Gauss map} is a continuous choice of a unit
normal $\nn:\Sigma\to \SS^2\subset \RR^3$. Observe that there are
two choices  $\nn$ and $-\nn$ corresponding to a choice of
orientation of $\Sigma$. If $\Sigma$ is minimal, then the Gauss
map is an (anti) conformal map since the eigenvalues of the
Weingarten map are $\kappa_1$ and $\kappa_2=-\kappa_1$. Moreover,
for a minimal surface
\begin{equation}        \label{e:gmconf}
        |A|^2=\kappa_1^2+\kappa_2^2=-2\,\kappa_1\,\kappa_2=-2\,K_{\Sigma} \, ,
\end{equation}
where $K_{\Sigma}$ is the Gauss curvature.  It follows that the
area of the image of the Gauss map is a multiple of the total
curvature.

\subsection{Minimal graphs}

Suppose that $u:\Omega\subset \RR^2 \to \RR$ is a $C^2$ function.
The graph of   $u$
\begin{equation}
        \Graph_u= \{(x,y,u(x,y)) \mid (x,y)\in \Omega\}\, .
\end{equation}
has area
\begin{align}
        \Area (\Graph_u)&=\int_{\Omega} |(1,0,u_x)\times (0,1,u_y)|\\
        &=\int_{\Omega} \sqrt{1+u_x^2+u_y^2}=\int_{\Omega}\sqrt{1+|\nabla u|^2}
        \, ,\notag
\end{align}
and the (upward pointing) unit normal is
\begin{equation}                \label{e:bb1}
        \nn=\frac{(1,0,u_x)\times (0,1,u_y)}{|(1,0,u_x)\times (0,1,u_y)|}
        =\frac{(-u_x,-u_y,1)}{\sqrt{1+|\nabla u|^2}}\, .
\end{equation}
Therefore for the graphs $\Graph_{u+t\eta}$ where $\eta|\partial
\Omega=0$ we get that
\begin{equation}
        \Area (\Graph_{u+t\eta})
        =\int_{\Omega}\sqrt{1+|\nabla u+t\,\nabla \eta|^2} \, ,
\end{equation}
hence
\begin{align}
       \left. \frac{d}{dt} \right|_{t=0} \Area (\Graph_{u+t\eta})
        &=\int_{\Omega} \frac{\langle \nabla u \, , \nabla \eta \rangle}{\sqrt{1+|\nabla u|^2}}\\
&=-\int_{\Omega}\eta\, \dv\, \left( \frac{\nabla
u}{\sqrt{1+|\nabla u|^2}} \right)\, .\notag
\end{align}
It follows that the graph of $u$ is a critical point for the area
functional if and only if $u$ satisfies the divergence form
equation
\begin{equation}  \label{e:mineq}
        \dv \left( \frac{\nabla u}{\sqrt{1+|\nabla u|^2}} \right)=0  \, .
\end{equation}
We will refer to \eqr{e:mineq} as the {\it minimal surface
equation}.

Next we want to show that the graph of a function on $\Omega$
satisfying the minimal surface equation  is not just a critical
point for the area functional but is actually area-minimizing
amongst surfaces in the cylinder $\Omega \times \RR \subset
\RR^3$.  To show this, first extend first the unit normal $\nn$ of
the graph in \eqr{e:bb1} to a vector field, still denoted by
$\nn$, on the entire cylinder $\Omega \times \RR$ by setting
\begin{equation}
    \nn (x,y,z) = \nn
(x,y, u(x,y)) \, .
\end{equation}
 Let $\omega$ be the two-form on
$\Omega \times \RR$ given by the condition that for $X,\, Y\in
\RR^3$
\begin{equation}
\omega (X,Y)=\text{det} (X,Y,\nn)\, .
\end{equation}
An easy calculation shows that
\begin{equation}
d \omega=\frac{\partial }{\partial x} \left(
\frac{-u_x}{\sqrt{1+|\nabla u|^2}} \right) +\frac{\partial
}{\partial y} \left( \frac{-u_y}{\sqrt{1+|\nabla u|^2}} \right)
=0\, ,
\end{equation}
since $u$ satisfies the minimal surface equation. In sum, the
form $\omega$ is closed and, given any
  $X$ and $Y$ at a point $(x,y,z)$,
\begin{equation}        \label{e:calibration1}
        | \omega (X,Y) | \leq | X \times Y| \, ,
\end{equation}
where  equality holds if and only if
\begin{equation}        \label{e:calibration2}
        X, Y \subset T_{(x,y, u(x,y))}
        \Graph_u \, .
\end{equation}
  Such a form $\omega$ is called a {\it calibration}.
  From this, we have that if $\Sigma \subset \Omega \times \RR$
is any other surface with $\partial \Sigma=\partial\, \Graph_u$,
then by Stokes' theorem since $\omega$ is closed,
\begin{equation}  \label{e:areamin}
        \Area (\Graph_u) = \int_{\Graph_u} \omega =
                \int_{\Sigma} \omega \leq
        \Area (\Sigma) \,  .
\end{equation}
This shows that $\Graph_u$ is area-minimizing among all surfaces
in the cylinder and with the same boundary.

If the domain $\Omega$ is convex, the minimal graph is absolutely
area-minimizing. To see this, observe first that if $\Omega$ is
convex, then so is $\Omega \times \RR $, and hence the nearest
point projection $P: \RR^3 \to \Omega \times \RR$ is a distance
nonincreasing Lipschitz map that is equal to the identity on
 $\Omega \times \RR$.
If $\Sigma \subset \RR^3$ is any other surface with $\partial
\Sigma=\partial\, \Graph_u$, then $\Sigma' = P(\Sigma)$ has $\Area
(\Sigma') \leq  \Area (\Sigma)$.  Applying \eqr{e:areamin} to
$\Sigma'$, we see that $\Area   (\Graph_u) \leq \Area (\Sigma')$
and the claim follows.

If $\Omega\subset \RR^2$ contains a ball of radius $r$, then,
since $\partial B_r\cap \Graph_u$ divides $\partial B_r$ into two
components at least one of which has area at most  $ ( \Area
(\SS^2) / 2 ) \,r^2$, we get from \eqr{e:areamin} the crude
estimate
\begin{equation}                \label{e:area}
        \Area (B_r\cap \Graph_u)\leq \frac{\Area (\SS^2)}{2}\, r^2\, .
\end{equation}

Very similar calculations to the ones above show that if
$\Omega\subset \RR^{n-1}$ and $u:\Omega\to \RR$ is a $C^2$
function, then the graph of $u$ is a critical point for the area
functional if and only if $u$ satisfies \eqr{e:mineq}. Moreover,
as in \eqr{e:areamin}, the graph of $u$ is actually
area-minimizing. Consequently, as in \eqr{e:area}, if $\Omega$
contains a ball of radius $r$, then
\begin{equation}                \label{e:areaforn}
        \Vol (B_r\cap \Graph_u)\leq \frac{\Vol (\SS^{n-1})}{2}\, r^{n-1}\, .
\end{equation}

\subsection{The maximum principle}

The first variation formula, \eqr{e:frstvar}, showed that a smooth
submanifold is a critical point for area if and only if the mean
curvature vanishes.  We will next
  derive the weak form of the  first variation formula
which is the basic tool for working with ``weak solutions''
(typically, stationary varifolds). Let $X$ be a vector field on
$\RR^n$.  We can write the divergence $\dv_{\Sigma} \, X$ of $X$
on $\Sigma$ as
\begin{equation}  \label{e:o1.3.1}
        \dv_{\Sigma} \, X=\dv_{\Sigma}\, X^T + \dv_{\Sigma}\, X^N = \dv_{\Sigma}\, X^T
        + \langle X , \HH \rangle \, ,
\end{equation}
where $X^T$ and $X^N$ are the tangential and normal projections of
$X$.  In particular, we get that, for a minimal submanifold,
\begin{equation}  \label{e:o1.3.1a}
        \dv_{\Sigma} \, X=\dv_{\Sigma}\, X^T
       \, .
\end{equation}
Moreover, from \eqr{e:o1.3.1} and Stokes' theorem, we see that
$\Sigma$ is minimal if and only if for all vector fields $X$ with
compact support and vanishing on the boundary of $\Sigma$,
\begin{equation}   \label{e:o1.3.2}
        \int_{\Sigma}\dv_{\Sigma}\, X=0\, .
\end{equation}
The key point is that \eqr{e:o1.3.2} makes sense as long as we can
define the divergence of a vector field on $\Sigma$.

The two most common notions of weakly minimal submanifolds are
minimal currents and stationary varifolds.   A {\it $k$-varifold}
$\Sigma$ on $\RR^n$ is a Radon measure on $\RR^n \times G(k,n)$,
where
 $G(k,n)$ is the space of (unoriented) $k$-planes through the
origin in $\RR^n$ (so that $G(n-1,n)$ is projective $(n-1)$
space).   If $X$ is a $C^1$ vector field, then we define the
divergence with respect to $(x,\omega) \in \RR^n \times G(k,n)$ by
\begin{equation}        \label{e:fvv2}
        \dv_{\omega} X
          = \sum_{i=1}^k \langle E_i , \nabla_{E_i} X \rangle \, ,
\end{equation}
where $E_i$ is an orthonormal basis for $\omega$ and $\nabla$ is
the Euclidean derivative.  In particular, the equation
\eqr{e:o1.3.2} makes sense when $\Sigma$ is a varifold and we will
say that a varifold satisfying \eqr{e:o1.3.2} is {\it stationary}.

 As a consequence of
\eqr{e:o1.3.2},
 we have the following proposition:

\begin{proposition}
\label{p:haco} $\Sigma^k\subset \RR^n$ is minimal if and only if
the restrictions of the coordinate functions of $\RR^n$ to
$\Sigma$ are harmonic functions.
\end{proposition}

\begin{proof}
Let $\eta$ be a smooth function on $\Sigma$ with compact support
and $\eta |\partial \Sigma=0$, then
\begin{equation}
        \int_{\Sigma} \langle \nabla_{\Sigma} \eta , \,
        \nabla_{\Sigma} x_i \rangle = \int_{\Sigma} \langle  \nabla_{\Sigma}\eta \, ,
        e_i \rangle =  \int_{\Sigma} \dv_{\Sigma} (\eta\, e_i)\, .
\end{equation}
From this, the claim follows easily.
\end{proof}

Recall that if $\Xi\subset \RR^n$ is a compact subset, then the
smallest convex set containing $\Xi$ (the convex hull, $\Conv
(\Xi)$) is the intersection of all half--spaces containing $\Xi$.
The maximum principle forces a compact minimal submanifold to lie
in the convex hull of its boundary  (this is the  ``convex hull
property''):

\begin{proposition}
If $\Sigma^k\subset \RR^n$ is a compact minimal submanifold,
 then
$\Sigma\subset \Conv (\partial \Sigma)$.
\end{proposition}

\begin{proof}
A half--space $H \subset \RR^n$ can be written as
\begin{equation}
H = \{ x \in \RR^n \, | \, \langle x , e \rangle \leq a \} \, ,
\end{equation}
 for a vector $e \in \SS^{n-1}$ and constant $a \in \RR$.
  By Proposition \ref{p:haco},   the function
$u(x)=\langle e, x\rangle$ is harmonic on $\Sigma$ and hence
attains its maximum on $\partial \Sigma$ by the maximum principle.
\end{proof}

Another application of \eqr{e:o1.3.1a}, with a different choice of
vector field $X$, gives that for a $k$-dimensional minimal
submanifold $\Sigma$
\begin{equation}
        \Delta_{\Sigma} |x-x_0|^2=2 \, \dv_{\Sigma} (x-x_0)=2k \, .
\end{equation}
Later we will see that this formula plays a crucial role in the
monotonicity formula for minimal submanifolds.

\vskip2mm The argument in the proof of the convex hull property
can be rephrased as saying that as we translate a hyperplane
towards a minimal surface, the first point of contact must be on
the boundary.  When $\Sigma$ is a hypersurface, this is a special
case of the strong maximum principle for minimal surfaces (see
\cite{CM1} for a proof):

\begin{lemma}             \label{l:smp}
Let $\Omega \subset \RR^{n-1}$ be an open connected neighborhood
of the origin. If $u_1$, $u_2:\Omega\to \RR$ are solutions of the
minimal surface equation with $u_1\leq u_2$ and $u_1(0)=u_2(0)$,
then $u_1\equiv u_2$.
\end{lemma}

Since any smooth hypersurface is locally a graph over a
hyperplane, Lemma \ref{l:smp} gives a maximum principle for smooth
minimal hypersurfaces.  There have been several interesting
extensions of the maximum principle to the singular case.  For
example, L. Simon, \cite{Si4}, proved that the strong maximum
holds when both hypersurfaces are area-minimizing even in the
presence of singularities.
  B. Solomon and B. White,
\cite{SoWh},  showed that it holds when at least one of the
hypersurfaces is smooth (the other may be just a stationary
varifold).  Finally, T. Ilmanen, \cite{I}, showed that it holds
for stationary hypersurfaces as long as each singular set has zero
codimension two measure.

 \vskip2mm  Thus far, the examples of minimal submanifolds have all been smooth.
The simplest non-smooth example
 is given by a pair of planes intersecting transversely along a
 line.  To get an example that is not even immersed, one can take
 three half--planes meeting along a line with an
 angle of
 $2\pi/3$ between each adjacent pair.

\section{Monotonicity and the mean value inequality}
\label{s:4}

Monotonicity formulas and mean value inequalities play a
fundamental role in many areas of geometric analysis.  Before we
state and prove the monotonicity formula of volume for minimal
submanifolds, we will need to recall the coarea formula. This
formula asserts (see, for instance, \cite{Fe} for a proof) that if
$\Sigma$ is a manifold and $h:\Sigma\to \RR$ is a  Lipschitz
function on $\Sigma$, then for all locally integrable functions
$f$ on $\Sigma$ and $t\in \RR$
\begin{equation}                \label{e:thecoareaf}
        \int_{\{h\leq t\}} f\,|\nabla h|
                =\int_{-\infty}^t \int_{h=\tau} f\,d\tau\, .
\end{equation}
When we apply this formula below, $h$ will be Euclidean distance
to a fixed point $x_0$.

\begin{proposition}
\label{p:pmonot} Suppose that $\Sigma^k\subset \RR^n$ is a minimal
submanifold and $x_0\in \RR^n$; then for all $0<s<t$
\begin{equation}        \label{e:monot0}
        t^{-k}\, \Vol (B_t (x_0) \cap \Sigma) -s^{-k}\, \Vol (B_s (x_0) \cap \Sigma)
        =\int_{(B_t(x_0) \setminus B_s(x_0) )\cap \Sigma} \frac{|(x-x_0)^N|^2}{|x-x_0|^{k+2}}\, .
\end{equation}
\end{proposition}

\begin{proof}
Within this proof, we set $B_t = B_t (x_0)$.  Since $\Sigma$ is
minimal,
\begin{equation}
        \Delta_{\Sigma} |x-x_0|^2=2 \, \dv_{\Sigma} (x-x_0)=2k \, .
\end{equation}
By Stokes' theorem integrating this gives
\begin{equation}
2\, k\, \Vol (B_s\cap \Sigma)=\int_{B_s\cap \Sigma}
\Delta_{\Sigma} |x-x_0|^2 =2\int_{\partial B_s\cap \Sigma}
|(x-x_0)^T|\, .
\end{equation}
Using this and the coarea formula (i.e., \eqr{e:thecoareaf}), an
easy calculation gives
\begin{align}
        \frac{d}{ds}\left( s^{-k}\, \Vol (B_s\cap \Sigma)\right)
        &=-k\, s^{-k-1}\, \Vol (B_s\cap \Sigma)
        +s^{-k}\,\int_{\partial B_s\cap \Sigma}
                \frac{|x-x_0|}{|(x-x_0)^T|}\notag\\
        &=s^{-k-1}\int_{\partial B_s\cap \Sigma}
        \left( \frac{|x-x_0|^2}{|(x-x_0)^T|}-|(x-x_0)^T| \right)\\
        &=s^{-k-1}\int_{\partial B_s\cap \Sigma}
                \frac{|(x-x_0)^N|^2}{|(x-x_0)^T|}\, .
        \notag
\end{align}
Integrating and applying the coarea formula once more gives the
claim.
\end{proof}

Notice that $(x-x_0)^N$ vanishes precisely when $\Sigma$ is
conical about $x_0$, i.e., when $\Sigma$ is invariant under
dilations about $x_0$. As a corollary, we get the following:

\begin{corollary}             \label{c:cmon}
Suppose that $\Sigma^k\subset \RR^n$ is a minimal submanifold and
$x_0\in \RR^n$; then the function
\begin{equation}                \label{e:thetadef}
        \Theta_{x_0}(s)
        =\frac{\Vol (B_s (x_0) \cap \Sigma)}{\Vol (B_s \subset \RR^k) }\,
\end{equation}
is a nondecreasing function of $s$.  Moreover, $\Theta_{x_0}(s)$
is constant in $s$ if and only if $\Sigma$ is conical about $x_0$.
\end{corollary}

Of course, if $x_0$ is a smooth point of $\Sigma$, then
$\lim_{s\to 0} \Theta_{x_0}(s) = 1$.  We will later see that the
converse is also true; this will be a consequence of the Allard
regularity theorem.

The monotonicity of area is a very useful tool in the regularity
theory for minimal surfaces --- at least when there is some {\it a
priori} area bound.  For instance, this monotonicity and a
compactness argument allow one to reduce many regularity questions
to questions about minimal cones (this was a key observation of W.
Fleming in his work on the Bernstein problem; see Section
\ref{s:5}).  Similar monotonicity formulas have played key roles
in other geometric problems, including  harmonic maps, Yang--Mills
connections, J--holomorphic curves, and regularity of limit spaces
with a lower Ricci curvature bound.

\vskip2mm  Arguing as in Proposition \ref{p:pmonot}, we get a
weighted monotonicity:

\begin{proposition}
\label{p:meanvalue} If\hspace{2pt}  $\Sigma^k\subset \RR^n$ is a
minimal submanifold, $x_0\in \RR^n$, and $f$ is a twice
differentiable function on $\Sigma$, then
\begin{equation}
        t^{-k}\int_{B_t(x_0) \cap \Sigma} f-s^{-k}\int_{B_s(x_0) \cap \Sigma} f
\end{equation}
\begin{equation}
        =\int_{(B_t(x_0)\setminus B_s(x_0))\cap \Sigma} f\,
                \frac{|(x-x_0)^N|^2}{|x-x_0|^{k+2}}
        +\frac{1}{2} \int_{s}^t \tau^{-k-1}
        \int_{B_{\tau}(x_0)\cap \Sigma} (\tau^2-|x-x_0|^2)\,
                \Delta_{\Sigma} f\,d\tau
        \, . \notag
\end{equation}
\end{proposition}

We get immediately the following mean value inequality for the
special case of non--negative subharmonic functions:

\begin{corollary}
Suppose that $\Sigma^k\subset \RR^n$ is a minimal submanifold,
$x_0\in \RR^n$, and $f$ is a non--negative subharmonic function on
$\Sigma$;
 then
\begin{equation}                \label{e:defavg}
        s^{-k}  \int_{B_s(x_0) \cap \Sigma} f
\end{equation}
is a nondecreasing function of $s$. In particular, if $x_0\in
\Sigma$, then for all $s>0$
\begin{equation}
        f(x_0)\leq
        \frac{\int_{B_s (x_0) \cap \Sigma} f}
                {\Vol \, (B_s \subset \RR^k)} \, .
\end{equation}
\end{corollary}

\section{Rado's theorem}

One of the most basic questions is what does the boundary
$\partial \Sigma$ tell us about a compact minimal submanifold
$\Sigma$?  We have already seen that $\Sigma$ must lie in the
convex hull of $\partial \Sigma$, but there are many other
theorems of this nature.  One of the first is a beautiful result
of T. Rado which says that if $\partial \Sigma$ is a graph over
the boundary of a convex set in $\RR^2$, then $\Sigma$ is also
graph (and hence embedded). The proof of this uses basic
properties of nodal lines for harmonic functions.

\begin{theorem}            \label{t:rado}
Suppose that $\Omega\subset \RR^2$ is a convex subset and
$\sigma\subset \RR^3$ is a simple closed curve which is graphical
over $\partial \Omega$.  Then any minimal disk $\Sigma \subset
\RR^3$ with $\partial \Sigma =\sigma$ must be graphical over
$\Omega$ and hence unique by the maximum principle.
\end{theorem}

\begin{proof}
(Sketch.) The proof is by contradiction, so suppose that $\Sigma$
is such a minimal disk  and $x \in \Sigma$ is a point where
  the tangent plane to
$\Sigma$  is vertical.  Consequently, there exists $(a,b)\ne
(0,0)$ such that
\begin{equation}        \label{e:rado1}
        \nabla_{\Sigma} (a\,x_1+b\,x_2)(x) = 0 \, .
\end{equation}
By Proposition \ref{p:haco}, $a\,x_1+b\,x_2$ is harmonic on
$\Sigma$ (since it is a linear combination of coordinate
functions).  The local structure of nodal sets of harmonic
functions (see, e.g., \cite{CM1}) then gives that the level set
\begin{equation}                \label{e:rado2}
        \{ y \in \Sigma \, | \,  (a\,x_1+b\,x_2) (y)
        = (a\,x_1+b\,x_2) (x) \}
\end{equation}
 has a singularity
at $x$ where at least four different curves meet. If two of these
nodal curves were to meet again, then there would be a closed
nodal curve which must bound a disk (since $\Sigma$ is a disk).
By the maximum principle, $a\,x_1 + b\, x_2$ would   have to
 be constant on this disk
and hence constant on $\Sigma$ by unique continuation. This would
imply that $\sigma = \partial \Sigma$ is contained in the plane
given by \eqr{e:rado2}.  Since this is impossible, we conclude
that all of these curves go to the boundary without intersecting
again.

In other words, the plane in $\RR^3$ given by \eqr{e:rado2}
intersects $\sigma$ in at least four points.  However, since
 $\Omega \subset \RR^2$ is convex, $\partial \Omega$ intersects
 the line given
by \eqr{e:rado2} in exactly two points.  Finally, since $\sigma$
is graphical  over $\partial \Omega$, $\sigma$ intersects the
plane in $\RR^3$ given by \eqr{e:rado2} in exactly two points,
which gives the desired contradiction.
\end{proof}

\section{The theorems of Bernstein and Bers}        \label{s:5}

 A classical theorem of S. Bernstein from 1916 says that
entire (i.e., defined over all of $\RR^2$) minimal graphs are
planes. This remarkable theorem of Bernstein was one of the first
illustrations of the fact that the solutions to a nonlinear PDE,
like the minimal surface equation, can behave quite differently
from solutions to a linear equation.

\begin{theorem} \cite{Be}  \label{t:bern} If
$u:\RR^{2}\to \RR$ is an entire solution to the minimal surface
equation, then $u$ is an affine function.
\end{theorem}

\begin{proof}
(Sketch.)  We will show that the curvature of the graph vanishes
identically; this implies that the unit normal is constant and,
hence, the graph must be a plane.  The proof follows by combining
two facts.  First, the area estimate for graphs \eqr{e:area} gives
\begin{equation}    \label{e:qaggg}
        \Area (B_r\cap \Graph_u)\leq 2 \, \pi \, r^2\, .
\end{equation}
This quadratic area growth allows one to construct a sequence of
non-negative logarithmic cutoff functions $\phi_j$ defined on the
graph with $\phi_j \to 1$ everywhere and
\begin{equation}    \label{e:tozero}
   \lim_{j\to \infty} \,  \int_{\Graph_u} |\nabla \phi_j|^2 = 0 \, .
\end{equation}
Moreover, since graphs are area-minimizing, they must be stable.
We can therefore use $\phi_j$ in the stability inequality
\eqr{e:stabin} to get
\begin{equation}    \label{e:stabiii}
    \int_{\Graph_u} \phi_j^2 \, |A|^2 \leq \int_{\Graph_u} |\nabla \phi_j|^2  \, .
\end{equation}
Combining these gives that $|A|^2$ is zero, as desired.
\end{proof}

 Rather surprisingly,
this result very much depended on the dimension.  The combined
efforts of E. De Giorgi \cite{DG}, F. J. Almgren, Jr. \cite{Am1},
and J. Simons \cite{Sim} finally gave:

\begin{theorem}  If
$u:\RR^{n-1}\to \RR$ is an entire solution to the minimal surface
equation and $n\leq 8$, then $u$ is an affine function.
\end{theorem}

However, in 1969 E. Bombieri, De Giorgi, and E. Giusti \cite{BDGG}
constructed entire non--affine solutions to the minimal surface
equation on $\RR^8$ and an area--minimizing singular cone in
$\RR^8$. In fact, they showed that for $m \geq 4$ the cones
\begin{equation}        \label{e:cones1}
        C_m = \{ (x_1 , \dots , x_{2m} ) \mid  x_1^2 + \cdots + x_m^2 =
x_{m+1}^2 + \cdots + x_{2m}^2  \} \subset \RR^{2m} \,
\end{equation}
are area--minimizing (and obviously singular at the origin).

\vskip2mm
 In contrast to the entire case, exterior solutions of the minimal graph
 equation, i.e., solutions on $\RR^2 \setminus B_1$, are much more
 plentiful. In this case,   L. Bers proved that $\nabla u$
actually has an asymptotic limit:

\begin{theorem} \cite{Ber} \label{t:bers}
If $u$ is a $C^2$ solution to the minimal surface equation on
$\RR^2 \setminus B_1$, then $\nabla u$ has a limit at infinity
(i.e., there is an asymptotic tangent plane).
\end{theorem}

\begin{proof}
(Sketch.) To get a rough idea of why Bers' theorem should hold,
recall that the Gauss map $\nn$ gives a holomorphic map from the
graph of $u$ to the upper hemisphere.  Hence, composing $\nn$ with
stereographic projection gives a bounded holomorphic function on
the graph of $u$.

The rest of the argument is to show that the point at infinity is
a removable singularity for this bounded holomorphic function.
Since the graph of $u$ is a topological annulus, it must be
conformal to an annulus of the form $\{ a < |z| \leq 1 \}$ in the
plane for some $a \geq 0$.  The point is to show that $a=0$.
 Arguing as in the proof of
Bernstein's theorem (cf. \eqr{e:qaggg}) shows that the graph of
$u$ has quadratic area growth and, hence, we can construct
logarithmic cutoff functions ``at infinity'' with arbitrarily
small energy.  It follows that we must have $a=0$ and hence $\nn$
extends smoothly across the puncture.  This gives a limiting value
for $\nn$ and, consequently, also a limiting value for $\nabla u$.
\end{proof}

 Bers' theorem was extended to higher
dimensions by L. Simon:

\begin{theorem} \cite{Si1} \label{t:berssi}
If $u$ is a $C^2$ solution to the minimal surface equation on
$\RR^n \setminus B_1$, then either \begin{itemize} \item $|\nabla
u|$ is bounded and $\nabla u$ has a limit at infinity. \item All
tangent cones at infinity are of the form $\Sigma \times \RR$
where $\Sigma$ is singular. \end{itemize}
\end{theorem}

Bernstein's theorem has had many other interesting
generalizations, some of which will be discussed later.

\section{Simons inequality}  \label{s:9}

In this section, we recall a very useful differential inequality
for the Laplacian of the norm squared of the second fundamental
form $A$ of a minimal hypersurface $\Sigma$ in $\RR^n$ and
illustrate its role in {\it a priori} estimates. This inequality,
originally due to J. Simons (see \cite{CM1} for a proof and
further discussion), is:

\begin{lemma}  \cite{Sim}   \label{l:simonsine}
If $\Sigma^{n-1}\subset \RR^n$ is a minimal hypersurface, then
\begin{equation}    \label{e:simtype}
\Delta_{\Sigma} \, |A|^2  = - 2\, |A|^4 + 2 |\nabla_{\Sigma} A|^2
\geq - 2 \, |A|^4 \, .
\end{equation}
\end{lemma}

An inequality of the type \eqr{e:simtype} on its own does not lead
to pointwise bounds on $|A|^2$ because of the nonlinearity.
However, it does lead to estimates if a   ``scale--invariant
energy'' is small.  For example, H. Choi and R. Schoen used
\eqr{e:simtype}  to prove:

\begin{theorem} \cite{CiSc}     \label{t:cisc}
There exists $\epsilon > 0$ so that if $0 \in \Sigma \subset B_r
(0)$ with $\partial \Sigma \subset
\partial B_r (0)$ is a minimal surface with
\begin{equation}
    \int |A|^2 \leq \epsilon \, ,
\end{equation}
 then
\begin{equation}
    |A|^2(0) \leq   r^{-2} \, .
\end{equation}
\end{theorem}

\section{Heinz's curvature estimate for graphs}

One of the key themes in minimal surface theory is the usefulness
of {\it a priori} estimates.
   A basic example is the curvature estimate of
E. Heinz for minimal graphs over a disk $D_{r_0}$ of radius $r_0$
in the plane.

\begin{theorem}   \cite{He}           \label{t:heinz}
If $D_{r_0} \subset \RR^2$ and
 $u:D_{r_0}\to \RR$ satisfies
the minimal surface equation, then for $\Sigma={\text{Graph}}_u$
and $0<\sigma\leq r_0$
\begin{equation}        \label{e:heinz}
        \sigma^2\,\sup_{D_{r_0-\sigma}}|A|^2\leq C\, .
\end{equation}
\end{theorem}

\begin{proof}
(Sketch.)  Observe first that it suffices to prove the estimate
for $\sigma = r_0$, i.e., to show that
\begin{equation}        \label{e:heinz2}
        |A|^2 (0,u(0)) \leq C \, r_0^{-2} \, .
\end{equation}
Recall that minimal graphs are automatically stable. As in the
proof of Theorem \ref{t:bern}, the area estimate for graphs
\eqr{e:area} allows us to use a logarithmic cutoff function in the
the stability inequality \eqr{e:stabin} to get that
\begin{equation}
    \int_{B_{r_1} \cap \Graph_u}  |A|^2 \leq \frac{ C }{ \log (r_0/r_1) } \, .
\end{equation}
Taking $r_0/ r_1$ sufficiently large, we can then apply Theorem
\ref{t:cisc} to get \eqr{e:heinz2}.
\end{proof}

 Heinz's estimate gives an effective version
of the Bernstein's theorem; namely, letting the radius $r_0$ go to
infinity in \eqr{e:heinz} implies that $|A|$ vanishes, thus giving
Bernstein's theorem.

\section{Embedded minimal disks with area bounds}

In the early nineteen--eighties R. Schoen and L. Simon extended
the theorem of Bernstein to complete simply connected embedded
minimal surfaces in $\RR^3$ with quadratic area growth; see
\cite{ScSi}. A surface $\Sigma$ is said to have {\it quadratic
area growth} if for all $r>0$, the area of the intersection of the
surface with the ball in $\RR^3$ of radius $r$ and center at the
origin is bounded by $C\, r^2$ for a fixed constant $C$
independent of $r$.

\begin{theorem}      \cite{ScSi}  \label{t:quasgm}
Let $0 \in \Sigma^{2} \subset B_{r_0} = B_{r_0}(x) \subset \RR^3$
be an embedded simply connected minimal surface with $\partial
\Sigma \subset \partial B_{r_0}$. If $\mu > 0$ and either
\begin{align}   \label{e:case1}
        \Area (\Sigma) &\leq \mu \, r_0^2 \, , {\text{ or }}
\\
        \int_{\Sigma} |A|^2 &\leq \mu \, ,
\end{align}
then for the connected component $\Sigma'$ of $B_{r_0 / 2}(x_0)
\cap \Sigma$ with $0 \in \Sigma'$ we have
\begin{equation}        \label{e:oqgmc}
        \sup_{\Sigma'}  |A|^2
                \leq C \, r_0^{-2}
\end{equation}
for some $C=C(\mu)$.
\end{theorem}

 In corollary $1.18$ in
\cite{CM4}, this was generalized to quadratic area growth for
{\it{intrinsic}} balls (this generalization played an important
role in analyzing the local structure of embedded minimal
surfaces).  We will use $\cB_{r}(p)$ to denote the {\it intrinsic
ball} of radius $r$ centered at a point $p$ in a surface $\Sigma$,
i.e., $\cB_{r}(p)$ is the set of points $q$ in $\Sigma$   that can
be connected to $p$ by a path $\gamma$ {\underline{in $\Sigma$}}
of length less than $r$.

\begin{theorem} \cite{CM4}
Given a constant $C_I$, there exists $C_P$ so that if
$\cB_{2r_0}\subset \Sigma\subset \RR^3$ is an embedded minimal
disk satisfying either
\begin{align}
        \Area (\cB_{2r_0}) &\leq C_I \, r_0^2  {\text{ or }}
\\
        \int_{\cB_{2r_0}} |A|^2 &\leq C_I \, ,
\end{align}
then
\begin{equation}
    \sup_{\cB_{s}} |A|^2 \leq C_P \, s^{-2} \,
.
\end{equation}
\end{theorem}

As an immediate consequence, letting $r_0 \to \infty$ gives
Bernstein-type theorems for embedded simply connected minimal
surfaces with either bounded density or finite total curvature.

The classical minimal surfaces known as Enneper's surface and the
catenoid show that neither ``embedded'' nor ``simply-connected''
can be removed.   Enneper's surface is a complete immersed minimal
disk but is not flat and is not embedded (see \eqr{e:enneper} for
the definition of Enneper's surface).  This example shows that
embeddedness is essential for these estimates. Similarly, the
catenoid  shows that simply-connected is essential. The catenoid
is the minimal surface in $\RR^3$ given by
\begin{equation}
    \{ (\cosh s\, \cos t,\cosh s\, \sin t,s) \, | \,
    s,t\in\RR \} \, ;
\end{equation}
 see figure \ref{f:f12}.

\section{Stable minimal surfaces}

It turns out that stable minimal surfaces have {\it a priori}
estimates.
 Since minimal graphs are stable,
the estimates for stable surfaces can be thought of as
generalizations of the earlier estimates for graphs.  These
estimates have been widely applied and are particularly useful
when combined with existence results for stable surfaces (such as
the solution of the Plateau problem).
  The starting point for these estimates is that, as we saw in \eqr{e:secvar}, stable
minimal surfaces satisfy
 the
stability inequality
\begin{equation}   \label{e:stabin2}
      \int |A|^2 \, \phi^2 \leq   \int |\nabla \phi|^2 \, .
\end{equation}

We will mention two such estimates.  The first is R. Schoen's
curvature estimate for stable surfaces:

\begin{theorem}   \cite{Sc1}  \label{t:stable2}
There exists a constant $C$ so that if $\Sigma  \subset  \RR^3$ is
an immersed stable minimal surface with trivial normal bundle and
$\cB_{r_0} \subset \Sigma \setminus \partial \Sigma$, then
\begin{equation}        \label{e:o11}
        \sup_{ \cB_{r_0-\sigma} }  |A|^2
                \leq C \, \sigma^{-2}\, .
\end{equation}
\end{theorem}

The second is an estimate for the area and total curvature of a
stable surface; for simplicity, we will state only the area
estimate:

\begin{theorem}   \cite{CM2}  \label{t:stable3}
If $\Sigma  \subset  \RR^3$ is an immersed stable minimal surface
with trivial normal bundle and $\cB_{r_0} \subset \Sigma \setminus
\partial \Sigma$, then
\begin{equation}        \label{e:abb}
       \Area \, (\cB_{r_0}) \leq 4 \pi \, r_0^2 / 3 \, .
\end{equation}
\end{theorem}

As mentioned,  we can use \eqr{e:abb} to bound the energy of a
 cutoff function in the stability inequality and, thus,
bound the total curvature of sub-balls.  Combining this with the
curvature estimate of Theorem \ref{t:cisc} gives Theorem
\ref{t:stable2}; see \cite{CM2}.   Note that the bound \eqr{e:abb}
is surprisingly sharp; even when $\Sigma$ is a plane, the area is
$\pi r_0^2$.

We will explain next why Theorem \ref{t:stable3} holds.
 The stability inequality can be used to get upper bounds for the
total curvature in terms of the area of a minimal surface. On the
other hand, we can use either the Gauss--Bonnet theorem or the
Jacobi equation to get the opposite bound. Combining these two
bounds will give the {\it a priori} bound on the area   of
intrinsic balls in a stable surface. More precisely, integrating
the Jacobi equation and using the Gauss equation $K_{\Sigma} =
-|A|^2/2$ gives
\begin{equation}    \label{e:jacga}
 4\,(\Area \, (\cB_{R}) - \pi \, R^2) = 2
\int_{0}^{R} \int_0^{t} \int_{\cB_{s}} |A|^2
    = \int_{\cB_{R}}|A|^2\,(R-r)^2 \, .
\end{equation}
The second equality uses two integrations by parts (i.e.,
$\int_0^R f(t) \, g''(t) \, dt$ with $f(t) = \int_0^t \int_{\cB_s}
|A|^2$ and $g(t) = (R-t)^2$). See  corollary $1.7$ of \cite{CM4}
for further details.

Using $\phi = R-r$ (where $r(x) = \dist_{\Sigma}(0,x)$) in the
stability inequality gives
\begin{equation}    \label{e:jacgab}
 4\,(\Area \, (\cB_{R}) - \pi \, R^2)
    = \int_{\cB_{R}}|A|^2\,(R-r)^2 \leq \int_{\cB_R} |\nabla (R-r)|^2 =
    \Area \, (\cB_{R}) \, .
\end{equation}
Consequently,
 $\Area \, (\cB_{R}) \leq 4\pi \, R^2 / 3$.

\section{Regularity theory}        \label{s:18}

In this section, we survey some of the key  ideas in classical
regularity theory, such as the role of monotonicity, scaling,
$\epsilon$-regularity theorems (such as W. Allard's theorem) and
tangent cone analysis (such as F. Almgren's refinement of H.
Federer's dimension reducing).

The starting point for all of this is the monotonicity of volume
for a minimal $k$--dimensional submanifold $\Sigma$. Namely,
Corollary \ref{c:cmon}  gives that the {\it density}
\begin{equation}                \label{e:thetadefvv}
        \Theta_{x_0}(s)
        =\frac{\Vol (B_s (x_0) \cap \Sigma)}{\Vol (B_s \subset \RR^k)
        }
\end{equation}
is a monotone non-decreasing function of $s$.  Consequently, we
can define the density $\Theta_{x_0}$ at the point $x_0$ to be the
limit as $s \to 0$ of $\Theta_{x_0}(s)$.  It  also follows easily
from monotonicity that the density is semi--continuous as a
function of $x_0$.

\vskip2mm \subsection{$\epsilon$--regularity and the singular set}

An {\it $\epsilon$--regularity theorem} is a theorem  giving that
a weak (or generalized) solution is actually smooth at a point if
a scale--invariant energy is small enough there. The standard
example is the Allard regularity theorem:
\begin{theorem} \cite{Al}   \label{t:allard}
There exists $\delta = \delta (k , n) > 0$ such that if $\Sigma
\subset \RR^{n}$ is a $k$--rectifiable stationary varifold (with
density at least one a.e.), $x_0 \in \Sigma$, and
\begin{equation}
    \Theta_{x_0} = \lim_{r\to 0} \frac{ \Vol \, (B_r (x_0) \cap \Sigma)}
        {\Vol \, ( B_r \subset \RR^k )}
        < 1 + \delta \, ,
\end{equation}
 then $\Sigma$ is smooth in a neighborhood of $x_0$.
\end{theorem}

Similarly, the small total curvature estimate of Theorem
\ref{t:cisc} may be thought of as an $\epsilon$-regularity
theorem; in this case, the scale--invariant energy is $\int
|A|^2$.

\vskip2mm
 As an application of the  $\epsilon$--regularity
theorem, Theorem \ref{t:allard},  we can define the singular set
$\cS$ of $\Sigma$ by
\begin{equation}
    \cS = \{ x \in \Sigma \, | \, \Theta_x \geq 1 + \delta \} \, .
\end{equation}
It follows immediately from the semi--continuity of the density
  that $\cS$ is closed.  In order to bound the
size of the singular set (e.g., the Hausdorff measure), one
combines the $\epsilon$--regularity with simple covering
arguments.

  This preliminary analysis of the singular set can be
refined by doing a so--called tangent cone analysis.

\subsection{Tangent cone analysis}

It is not hard to see that scaling  preserves the space of
  minimal
submanifolds of $\RR^n$.  Namely, if $\Sigma$ is minimal, then so
is
\begin{equation}
     \Sigma_{y,\lambda}  = \{ y + \lambda^{-1}\, (x-y) \, | \, x \in \Sigma \} \, .
\end{equation}
(To see this, simply note that this scaling multiplies the
principal curvatures by $\lambda$.) Suppose now that we fix the
point $y$ and take a sequence $\lambda_j \to 0$.  The monotonicity
formula   bounds the density of the rescaled solution, allowing us
to extract a convergent subsequence and limit.  This limit, which
is called a {\it tangent cone} at $y$,   achieves equality in the
monotonicity formula and, hence, must be homogeneous (i.e.,
invariant under dilations about $y$).

The usefulness of tangent cone analysis in
 regularity theory is based on two key facts.  For simplicity, we
 illustrate these when $\Sigma \subset \RR^n$ is an area
 minimizing hypersurface.   First, if any tangent cone at $y$ is a
 hyperplane $\RR^{n-1}$, then $\Sigma$ is smooth in a neighborhood
 of $y$.  This follows easily from the Allard regularity theorem since
 the density at $y$ of the tangent cone is the same as the density
 at $y$ of $\Sigma$.  The second key fact, known as ``dimension
 reducing,'' is due to F. Almgren, \cite{Am2},
  and is a refinement of an argument of
 H. Federer.  To state this, we first stratify the singular set $\cS$
 of $\Sigma$ into subsets
 \begin{equation}
    \cS_0 \subset \cS_1 \subset \cdots \subset \cS_{n-2} \, ,
 \end{equation}
 where we define $\cS_i$ to be the set of points $y\in \cS$ so
 that  any linear space contained in any tangent cone at $y$ has
 dimension at most $i$.  (Note that $\cS_{n-1} = \emptyset$ by Allard's
 Theorem.)  The dimension reducing argument then gives that
 \begin{equation}       \label{e:dimred}
    {\text{dim}} \, (\cS_i) \leq i \, ,
 \end{equation}
 where dimension   means the Hausdorff dimension.  In particular, the
 solution of the Bernstein problem then gives codimension $7$
 regularity of $\Sigma$, i.e., ${\text{dim}} \, (\cS) \leq n-8$.
 See lecture $2$ in \cite{Si3} for a proof of \eqr{e:dimred}.

\vskip2mm Tangent cones produced by scalings as above may very
well depend on the particular convergent subsequence.  In some
cases, one can prove uniqueness of the tangent cone and this is
often quite useful (see, for instance, section $3.4$ in \cite{Si3}
for one such application).

\part{Embedded minimal surfaces}

\stepcounter{section}

Thus far, the results for embedded minimal surfaces have assumed
some additional {\it a priori} bound, such as bounds on area or
total curvature, and the proofs break down without these {\it a
priori} bounds.  In this part, we will focus on recent results for
embedded minimal surfaces without {\it a priori} bounds.

\subsection{Multi--valued graphs}      \label{s:11}

 We have earlier studied minimal graphs. It is useful also to consider multi-valued
 minimal graphs.   Intuitively, an
(embedded) multi--valued graph is a surface such that over each
point of the annulus, the surface consists of $N$ graphs. To make
this notion precise, let  $\cP$ be the universal cover of the
punctured plane $\CC\setminus \{0\}$ with global polar coordinates
$(\rho, \theta)$ so $\rho>0$ and $\theta\in \RR$.  An {\it
$N$--valued graph} on  $\CC\setminus \{0\}$ is a single valued
graph of a function $u$ over $\cP$.
 For working purposes, we generally
think of the intuitive picture of a multi--sheeted surface in
$\RR^3$, and we identify the single--valued graph over the
universal cover with its multi--valued image in $\RR^3$.

We will also need the notion of an $N$-valued graph over the
annulus $\{ r < \rho \leq s \} \subset \CC$.  In this case, the
function is single-valued on
\begin{equation}
    \{(\rho,\theta)\,|\,r< \rho\leq s\, ,\, |\theta|\leq
    N\,\pi\} \, .
\end{equation}

\begin{figure}[htbp]
    \setlength{\captionindent}{20pt}
    \begin{minipage}[t]{0.5\textwidth}
    \centering\includegraphics{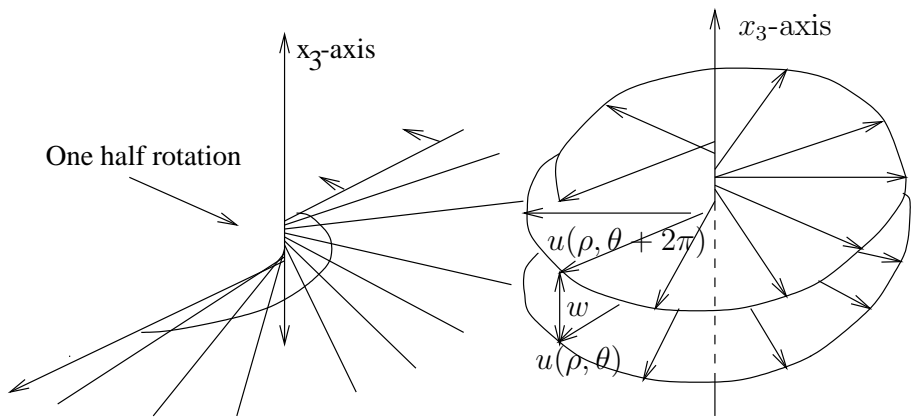}
    \caption{Multi--valued graphs.  The helicoid is obtained
    by gluing together two $\infty$--valued graphs along a line. First published in the Notices of the
    American Mathematical Society in 2003,
    published by the American Mathematical Society.}\label{f:f1}
    \end{minipage}%
    \begin{minipage}[t]{0.5\textwidth}
    \centering\input{pl8a.pstex_t}
    \caption{The separation $w$ grows/decays in $\rho$ at most sublinearly
for a multi--valued minimal graph; see \eqr{e:slgc}.}\label{f:f2}
    \end{minipage}
\end{figure}

The multi--valued graphs that we will consider will all be
embedded, which corresponds to a nonvanishing separation between
the sheets (or the floors).  Here the {\it separation} is the
function (see figure \ref{f:f2})
\begin{equation}
w(\rho,\theta)=u(\rho,\theta+2\pi)-u(\rho,\theta)\, .
\end{equation}

The helicoid is a minimal surface consisting of two multi-valued
graphs glued together along an axis.  The helicoid looks like a
``double spiral staircase'' (see \cite{CM10}) and is parametrized
by:

\vskip2mm \noindent {\bf{Example 2}}: (Helicoid; see figure
\ref{f:f1}).  The helicoid is the minimal surface in $\RR^3$ given
by the parametrization
\begin{equation}  \label{e:helicoid}
(s\cos t,s\sin t,t)\, ,\,\,\,\,\,\text{ where }s,\,t\in \RR\, .
\end{equation}

If $\Sigma$ is the helicoid, then $\Sigma\setminus
\{x_3-\text{axis}\}=\Sigma_1\cup \Sigma_2$,
 where $\Sigma_1$, $\Sigma_2$ are $\infty$--valued graphs on
$\CC\setminus \{0\}$. $\Sigma_1$ is the graph of the function
$u_1(\rho,\theta)=\theta$ and $\Sigma_2$ is the graph of the
function $u_2(\rho,\theta)=\theta+\pi$.  ($\Sigma_1$ is the subset
where $s>0$ in \eqr{e:helicoid} and $\Sigma_2$ the subset where
$s<0$.)  In either case the separation $w=2\,\pi$.  A {\it
multi--valued minimal graph} is a multi--valued graph of a
function $u$ satisfying the minimal surface equation.

Note that for an embedded multi--valued graph, the sign of $w$
determines whether the multi--valued graph spirals in a
left--handed or right--handed manner, in other words, whether
upwards motion corresponds to turning in a clockwise direction or
in a counterclockwise direction.

\subsection{The sublinear growth of the separation}

As we have seen, the separation is constant for the multi--valued
graphs coming from each half of the helicoid.  This can be viewed
as a type of Liouville Theorem reflecting the conformal properties
of an infinite--valued graph.  In Proposition II.2.12 of
\cite{CM3}, we proved a corresponding gradient estimate (i.e., a
pointwise estimate for $|\nabla \log |w||$) for a general
multi-valued graph. Integrating this gradient estimate gives  that
the separation grows sublinearly (see figure \ref{f:f2}):

\begin{proposition}     \label{c:slg}
\cite{CM3}
 Given $\alpha > 0$, there exists $N$ so that if
 $u$ satisfies the minimal surface equation
on \begin{equation} \label{e:thergg}
    \{ \e^{-N} \, r_1  \leq \rho \leq \e^{N} \,
    r_2 \, , \, - N \leq \theta \leq 2 \pi + N \} \, ,
\end{equation}
 $|\nabla u| \leq 1$, and has separation $w \ne 0$, then
\begin{equation}    \label{e:slgc}
    |w | (r_2,0)  \leq  |w| (r_1 , 0) \left( \frac{r_2}{r_1} \right)^{\alpha}   \, .
\end{equation}
\end{proposition}

\begin{proof}
(Sketch.)  As mentioned above, the inequality \eqr{e:slgc} follows
from integrating the gradient estimate
\begin{equation}    \label{e:graee}
    |\nabla \log |w|| (r,0) \leq \frac{\alpha}{r}
\end{equation}
  in the same way that the
Harnack inequality for positive harmonic functions follows from
integrating the gradient estimate.

To see why a gradient estimate like \eqr{e:graee} holds, observe
that $u (\cdot , \cdot)$ and its $2\pi$--rotation $u(\cdot , \cdot
+ 2\pi)$ are both solutions of the minimal surface equation and,
thus, the difference $w$ is (almost) a positive solution of the
linearized equation.  The linearized equation is itself a
perturbation of the Laplace equation and it is not difficult to
get an estimate
\begin{equation}    \label{e:graee1}
    |\nabla \log |w|| (r,0) \leq \frac{C}{r} \, ,
\end{equation}
for some constant $C$.

The point now is to show that if $N$ is large, then we can make
the constant $C$ in \eqr{e:graee1}  small.  For simplicity,
suppose that $w$ is actually a positive harmonic function   on the
region \eqr{e:thergg}.  Since harmonic functions are invariant
under conformal transformations, $w \circ \e^z$ is a positive
harmonic function on the ``rectangle''
\begin{equation}
    \{
(x+iy) \, | \, -N + \log r_1 < x < N  + \log r_2 {\text{ and }}
|y| < N \} \, .
\end{equation}
In the extreme case when $N=\infty$, the positive harmonic
function  $w \circ \e^z$ is defined on the entire plane and,
hence, is constant by the Liouville theorem. It is not hard to see
that applying the gradient estimate on this rectangle when $N$ is
large, and then translating this back to the original function
$w$, gives \eqr{e:graee}.
\end{proof}

\section{Embedded minimal disks}        \label{s:19}

There are two classical local models for embedded minimal disks
(by an {\it
  embedded disk} we mean a smooth injective map from the closed unit ball in
$\RR^2$ into $\RR^3$).  One model is the plane (or, more
generally, a minimal graph) and the other is a piece of a
helicoid.

The  helicoid  was discovered by Meusnier in 1776.  Meusnier had
been a student of Monge.  He also discovered that the surface now
known as the catenoid is minimal in the sense of Lagrange, and he
was the first to characterize a minimal surface as a surface with
vanishing mean curvature. Unlike the helicoid, the catenoid is not
topologically a plane but rather a cylinder.

\vskip2mm   It turns out that these two classical examples (graphs
and helicoids) completely capture the local structure of an
embedded minimal disk.  This is made concrete in the compactness
theorem for embedded minimal disks, Theorem \ref{t:t0.1} below.

  To avoid tedious dependence of various quantities we state
Theorem \ref{t:t0.1} not for a single embedded minimal disk with
sufficiently large curvature at a given point but instead for a
sequence of such disks where the curvatures are blowing up.
Theorem \ref{t:t0.1} says that a sequence of embedded minimal
disks mimics the following behavior of a sequence of rescaled
helicoids:

\begin{quotation}
Consider the sequence $\Sigma_i = a_i \, \Sigma$ of rescaled
helicoids where $a_i \to 0$. (That is, rescale $\RR^3$ by $a_i$,
so points that used to be distance $d$ apart will in the rescaled
$\RR^3$ be distance $a_i\,d$ apart.)  The curvatures of this
sequence of rescaled helicoids are blowing up along the vertical
axis. The sequence converges (away from the vertical axis) to a
foliation by flat parallel planes. The singular set $\cS$ (the
axis) then consists of removable singularities.
\end{quotation}

Let now $\Sigma_i\subset B_{2R}$ be a sequence of embedded minimal
 disks with
$\partial \Sigma_i\subset \partial B_{2R}$.  Clearly (after
possibly going to a subsequence) either (A) or (B) occurs:
\begin{enumerate}
\item[(A)]
 $\sup_{B_{R}\cap\Sigma_i}|A|^2\leq C<\infty$ for some constant $C$.
 \item[(B)]
 $\sup_{B_{R}\cap\Sigma_i}|A|^2\to \infty$.
\end{enumerate}
 In (A) (by a standard argument; see, e.g., lemma $2.2$ in \cite{CM1}) the intrinsic ball
$\cB_s(y_i)$ is a graph for all $y_i\in B_{R}\cap \Sigma_i$, where
$s$ depends only on $C$. Thus the main case is (B) which is the
subject of the next theorem.

\vskip6mm Using the notion of multi--valued graphs, the {\it
lamination theorem}  (the main theorem of \cite{CM6}), can now be
stated:

\begin{figure}[htbp]
    \setlength{\captionindent}{20pt}
    %\begin{minipage}[t]{0.5\textwidth}
    \centering\input{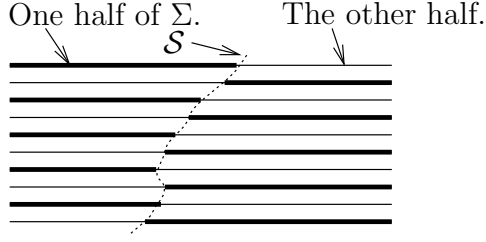}
    \caption{Theorem \ref{t:t0.1} -- the singular set, $\cS$, and
the two multi--valued graphs. First published in the Notices of
the American Mathematical Society in 2003,
    published by the American Mathematical Society.}\label{f:f3}
    %\end{minipage}
\end{figure}

\begin{theorem} \label{t:t0.1}
(Theorem 0.1 in \cite{CM6}).  See figure \ref{f:f3}. Let $\Sigma_i
\subset B_{R_i}=B_{R_i}(0)\subset \RR^3$ be a sequence of embedded
minimal disks with $\partial \Sigma_i\subset \partial B_{R_i}$
where $R_i\to \infty$. If
\begin{equation}
    \sup_{B_1\cap
\Sigma_i}|A|^2\to \infty \, ,
\end{equation}
 then there exists a
subsequence, $\Sigma_j$, and a Lipschitz curve $\cS:\RR\to \RR^3$
such that after a rotation of $\RR^3$:
\begin{enumerate}
\item
 $x_3(\cS(t))=t$.  (That is, $\cS$ is a
graph over the $x_3$--axis.)
\item  Each $\Sigma_j$
consists of exactly two multi--valued graphs away from $\cS$
(which spiral together). \item  For each $1>\alpha>0$,
$\Sigma_j\setminus \cS$ converges in the $C^{\alpha}$--topology to
the foliation, $\cF=\{x_3=t\}_t$, of $\RR^3$. \item
$\sup_{B_{r}(\cS (t))\cap \Sigma_j}|A|^2\to\infty$ for all $r>0$,
$t\in \RR$.  (The curvatures blow up along $\cS$.)
\end{enumerate}
\end{theorem}

In (2) and (3), the statements that $\Sigma_j\setminus \cS$ are
multi--valued graphs and converges to $\cF$ means that for each
compact subset $K\subset \RR^3\setminus \cS$ and $j$ sufficiently
large $K\cap \Sigma_j$ consists of multi--valued graphs over (part
of) $\{x_3=0\}$ and $K\cap \Sigma_j\to K\cap \cF$ in the sense of
graphs.

 A key point in  Theorem \ref{t:t0.1} is that there is no
useful monotonicity formula or natural {\it a priori} bound.  The
main tools for overcoming these difficulties are a
``classification of singularities'' which describes a neighborhood
of points of large curvature  and our one-sided curvature estimate
(Theorem \ref{t:t2} below).

\begin{figure}[htbp]
    \setlength{\captionindent}{4pt}
    %\begin{minipage}[t]{0.5\textwidth}
    \centering\input{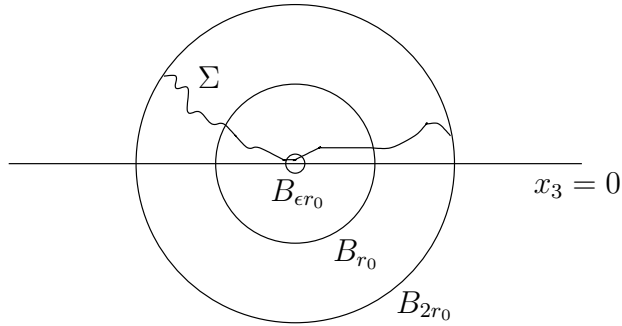}
    \caption{Theorem \ref{t:t2} -- the one--sided curvature estimate for an
embedded minimal disk $\Sigma$ in a half--space with $\partial
\Sigma\subset \partial B_{2r_0}$:  The components of $B_{r_0}\cap
\Sigma$ intersecting $B_{\epsilon r_0}$ are graphs. First
published in the Notices of the American Mathematical Society in
2003,
    published by the American Mathematical Society.}\label{f:f11}
    %\end{minipage}
\end{figure}

  The   one-sided curvature
estimate says roughly that if an embedded minimal disk lies in a
half--space above a plane and comes close to the plane, then it is
a graph over the plane. Precisely, this is the following theorem:

\begin{theorem}  \label{t:t2}
(Theorem 0.2 in \cite{CM6}). See figure \ref{f:f11}. There exists
$\epsilon>0$, so that if
$$\Sigma \subset B_{2r_0} \cap \{x_3>0\} \subset \RR^3$$ is an
embedded minimal disk with $\partial \Sigma\subset \partial B_{2
\, r_0}$, then for all components $\Sigma'$ of $B_{r_0} \cap
\Sigma$ which intersect $B_{\epsilon r_0}$ we have
\begin{equation}        \label{e:graph}
\sup_{\Sigma'} |A_{\Sigma}|^2 \leq r_0^{-2} \, .
\end{equation}
\end{theorem}

Using the minimal surface equation and that $\Sigma'$ has points
close to a plane, it is not hard to see that, for $\epsilon>0$
sufficiently small, \eqr{e:graph} is equivalent to the statement
that $\Sigma'$ is a graph over the plane $\{x_3=0\}$.

  We will often refer to
Theorem \ref{t:t2} as {\it the one--sided curvature estimate}.

\begin{figure}[htbp]
    \setlength{\captionindent}{20pt}
    \begin{minipage}[t]{0.5\textwidth}
    \centering\input{pl2a.pstex_t}
    \caption{The catenoid given by revolving $x_1= \cosh x_3$
around the $x_3$--axis. First published in the Notices of the
American Mathematical Society in 2003,
    published by the American Mathematical Society.}  \label{f:f12}
    \end{minipage}\begin{minipage}[t]{0.5\textwidth}
    \centering\input{unot7.pstex_t}
    \caption{Rescaling the catenoid shows that simply connected
(and embedded) is
    needed in the one--sided curvature estimate.
    First published in the Notices of the American Mathematical Society in 2003,
    published by the American Mathematical Society.}  \label{f:f13}
    \end{minipage}%
\end{figure}

Note that the assumption in Theorem \ref{t:t2} that $\Sigma$ is
simply connected is crucial as can be seen from the example of a
rescaled catenoid. Recall that the catenoid is the minimal surface
in $\RR^3$ given by
\begin{equation}
    (\cosh s\, \cos t,\cosh s\, \sin t,s)
\end{equation}
where $s,t\in\RR$; see figure \ref{f:f12}. Under rescalings this
converges (with multiplicity two) to the flat plane; see figure
\ref{f:f13}. Likewise, by considering the universal cover of the
catenoid, one sees that embedded, and not just immersed, is needed
in Theorem \ref{t:t2}.

As an almost immediate consequence of Theorem \ref{t:t2} and a
simple barrier argument we get that if in a ball two embedded
minimal disks come close to each other near the center of the ball
then each of the disks are graphs.  Precisely, this is the
following:

\begin{figure}[htbp]
    \setlength{\captionindent}{20pt}
    %\begin{minipage}[t]{0.5\textwidth}
    \centering\input{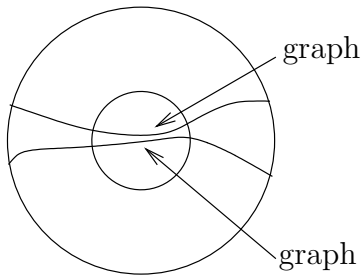}
    \caption{Corollary \ref{c:barrier}:  Two sufficiently close
    components of an embedded minimal disk must each be a graph.}\label{f:f14}
    %\end{minipage}
\end{figure}

\begin{corollary}         \label{c:barrier}
(Corollary 0.4 in \cite{CM6}). See figure \ref{f:f14}. There exist
$c > 1$, $\epsilon >0$ so that the following holds:

 \noindent Let
$\Sigma_1$ and $ \Sigma_2 \subset B_{cr_0} \subset \RR^3$ be
disjoint embedded minimal surfaces with $\partial \Sigma_i \subset
\partial B_{cr_0}$ and $B_{\epsilon \, r_0 } \cap \Sigma_i \ne
\emptyset$. If $\Sigma_1 $ is a disk,
 then for all components $\Sigma_1'$ of
$B_{r_0} \cap \Sigma_1$ which intersect $B_{\epsilon \, r_0}$
\begin{equation}        \label{e:onece}
    \sup_{\Sigma_1'}   |A|^2
        \leq  r_0^{-2}  \, .
\end{equation}
\end{corollary}

\section{Fixed genus}

Following our results on embedded minimal disks described in the
previous section, we proved
 two main structure theorems for {\it{non-simply connected}}
embedded minimal surfaces of any given fixed genus in \cite{CM7}.

The first of these asserts that any such surface
{\underline{without}} small necks can be obtained by gluing
together two oppositely--oriented double spiral staircases; see
figure \ref{f:01a}.

The second gives a pair of pants decomposition of any such surface
when there {\underline{are}} small necks, cutting the surface
along a collection of short curves; see figure \ref{f:01b}. After
the cutting, we are left with graphical pieces that are defined
over a disk with either one or two sub--disks removed (a
topological disk with two sub--disks removed is called a pair of
pants).

Both of these structures occur as different extremes in the
two-parameter family of minimal surfaces known as the Riemann
examples.

\begin{figure}[htbp]
    \center{The two main structure theorems for non-simply connected surfaces:}
        \vskip2mm
    \setlength{\captionindent}{20pt}
    \begin{minipage}[t]{0.5\textwidth}
    \centering\input{gen01a.pstex_t}
    \caption{Absence of necks:  The surface can be obtained by gluing together two
oppositely--oriented double spiral staircases.}
    \label{f:01a}
    \end{minipage}\begin{minipage}[t]{0.5\textwidth}
    \centering\input{gen01b.pstex_t}
    \caption{Presence of necks:  The surface can be decomposed into a collection of pair of pants
    by cutting along short curves.}
    \label{f:01b}
    \end{minipage}
\end{figure}

\subsection{Uniformly locally simply connected}

Sequences of surfaces which are not simply connected are, after
passing to a subsequence, naturally divided into two separate
cases depending on whether or not the topology is concentrating at
points. To distinguish between these cases, we will say that a
sequence of surfaces $\Sigma_i^2\subset \RR^3$ is {\it{uniformly
locally simply connected}} (or ULSC)  if for each compact subset
$K$ of $\RR^3$, there exists a constant $r_0 > 0$ (depending on
$K$) so that for every $x \in K$, all $r \leq r_0$, and every
surface $\Sigma_i$
\begin{equation}    \label{e:ulsc2}
 {\text{each connected component of }} B_{r}(x) \cap \Sigma_i {\text{ is
 a disk.}}
\end{equation}
  For instance, a sequence of rescaled catenoids
   where the necks shrink to zero is not ULSC, whereas a
sequence of rescaled helicoids is.

\begin{remark}
If each component of the  intersection of a minimal surface with a
ball of radius $r_0$ is a disk, then so are the intersections with
all sub--balls by  the convex hull property (see, e.g., lemma C.1
in \cite{CM6}). Therefore, it would be enough that \eqr{e:ulsc2}
holds for $r=r_0$.
\end{remark}

We will next describe briefly the case of ULSC sequences.  See
\cite{CM7} for more on this and for the general case of fixed
genus.

We will assume here that the surfaces are not disks  (the case of
disks was dealt with in the previous subsection). In particular,
we will assume that for each $i$, there exists some $y_i \in
\RR^3$ and $s_i>0$ so that
\begin{equation}    \label{e:notulsc}
{\text{ some component of }} B_{s_i}(y_i) \cap \Sigma_i {\text{
 is not a disk.}}
\end{equation}

\begin{figure}[htbp]
    \setlength{\captionindent}{20pt}
    \centering\input{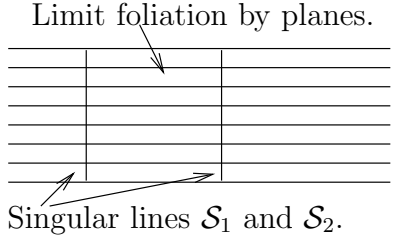}
    \caption{Theorem \ref{t:t5.1}: Limits of sequences of non-simply
connected, yet ULSC, surfaces with
    curvature blowing up.  The singular set consists of two lines
$\cS_1$ and $\cS_2$ and the limit  is a
    foliation by flat parallel planes.}
    \label{f:0}
\end{figure}

Loosely speaking, the next result  shows that when the sequence is
ULSC (but not simply connected), a subsequence converges to a
foliation by parallel planes away from two lines $\cS_1$ and
$\cS_2$; see figure \ref{f:0}. The lines $\cS_1$ and $\cS_2$ are
disjoint and orthogonal to the leaves of the foliation and the two
lines are precisely the points where the curvature is blowing up.
This is similar to the case of disks, except that we get two
singular curves for non-disks as opposed to just one singular
curve for disks.

\begin{theorem} \label{t:t5.1}  \cite{CM7}
 Let $\Sigma_i \subset B_{R_i}=B_{R_i}(0)\subset \RR^3$
be a sequence of compact embedded minimal surfaces with fixed
genus and with $\partial \Sigma_i\subset \partial B_{R_i}$ where
$R_i\to \infty$.  Suppose that each $\Sigma_i$ is ULSC and
satisfies
 \eqr{e:notulsc}  with
$s_i=R > 1$ and $y_i=0$. If
\begin{equation}
\sup_{B_1\cap \Sigma_i}|A|^2\to \infty \, ,
 \end{equation}
   then there exists a subsequence $\Sigma_j$, two  disjoint
  parallel
lines $\cS_1$ and $\cS_2$, and a rotation of $\RR^3$ so that:
 \begin{enumerate}
\item[(A)] For each $1>\alpha>0$, $\Sigma_j\setminus (\cS_1 \cup \cS_2)$
converges in the $C^{\alpha}$-topology to the foliation $\{ x_3 =
t \}$ by parallel planes.
\item[(B)]  $\sup_{B_{r}(x)\cap \Sigma_j}|A|^2\to\infty$ as $j \to \infty$ for
all $r>0$ and $x\in \cS_1 \cup \cS_2$.  (The curvatures blow up
along $\cS_1$ and $\cS_2$.)
\item[($C_{ulsc}$)]
Away from $\cS_1 \cup \cS_2$, each $\Sigma_j$ consists of exactly
two multi-valued graphs spiraling together.
 Near $\cS_1$ and $\cS_2$, the pair of multi-valued graphs form double
spiral staircases with opposite orientations at $\cS_1$ and
$\cS_2$. Thus, circling only $\cS_1$ or only $\cS_2$  results in
going either up or down, while a path circling both $\cS_1$ and
$\cS_2$ closes up (see figure \ref{f:mg2}).
\item[($D_{ulsc}$)] $\cS_1$ and $\cS_2$
are  orthogonal to the leaves of the foliation.
\end{enumerate}
\end{theorem}

We should point out that the multi-valued graphs in ($C_{ulsc}$)
are defined over the doubly-punctured plane; see figure
\ref{f:mg2}. The multi-valued graphs considered previously were
defined over a plane punctured at just one point.

\begin{figure}[htbp]
    \setlength{\captionindent}{20pt}
    \centering\input{gen1b.pstex_t}
    \caption{A multi-valued graph over the doubly-punctured plane.
    The spiral staircases near each puncture are oppositely--oriented.}
    \label{f:mg2}
\end{figure}

Despite the similarity  of Theorem \ref{t:t5.1}
 to the case of disks, it is worth noting that the results for
 disks do not alone give this theorem.  Namely, even though the
 ULSC sequence consists {\underline{locally}} of disks, the compactness
result for disks
 was in the {\underline{global}} case where the radii go to
 infinity.  One might wrongly think that Theorem \ref{t:t5.1} could be proven
 using the results for disks and a blow up argument. However, local
 examples constructed in \cite{CM11} show the difficulty with such an
 argument.

\subsection{Fixed genus}

In the previous section, we analyzed sequences of embedded minimal
surfaces with fixed genus where the topology does not concentrate
at any points (i.e., ULSC sequences).  We will now consider
general sequences where the topology is allowed to concentrate.
For simplicity, we will restrict to the case of genus zero, i.e.,
to planar domains (see \cite{CM7} for the general case).

One way of locally distinguishing sequences where the topology
does not concentrate from sequences where it does comes from
analyzing the singular set. The singular set $\cS$ is defined to
be the set of points where the curvature is blowing up. That is, a
point $y$ in $\RR^3$ is in $\cS$ for a sequence $\Sigma_i$ if
\begin{equation}
    \sup_{B_{r}(y)\cap
\Sigma_i}|A|^2\to\infty {\text{ as $i \to \infty$ for all $r>0$}}
.
\end{equation}
For embedded minimal surfaces, $\cS$ consists of two types of
points. The first type is roughly modelled on rescaled helicoids
 and the second on
rescaled catenoids:
\begin{itemize}
\item
A point $y$ in $\RR^3$ is in $\cSu$ if the curvature for the
sequence $\Sigma_i$ blows up at $y$ and the sequence is ULSC in a
neighborhood of $y$.
\item
A point $y$ in $\RR^3$ is in $\cSt$ if the sequence is not ULSC in
any neighborhood of $y$. In this case, a sequence of closed
non-contractible curves $\gamma_i \subset \Sigma_i$ converges to
$y$.
\end{itemize}
The sets $\cSt$ and $\cSu$ are obviously disjoint and
  the curvature blows up at both, so   $\cSt \cup \cSu \subset \cS$.   An easy
argument  shows that, after passing to a subsequence, we can
assume that
\begin{equation}    \label{e:}
    \cS=  \cSt \cup \cSu \, .
\end{equation}
Note that $\cSt = \emptyset$ is equivalent to that the sequence is
ULSC as is the case for sequences of rescaled helicoids.  On the
other hand, $\cSu = \emptyset$ for sequences of rescaled
catenoids.

We will show that every sequence $\Sigma_i$ has a subsequence that
is either    ULSC or for which $\cSu$ is empty. This is the
 next ``no mixing'' theorem.  These two different cases give
 two very different structures.

\begin{theorem}     \label{c:main}
\cite{CM7}  If $\Sigma_i \subset B_{R_i}=B_{R_i}(0)\subset \RR^3$
is a sequence
 of compact embedded minimal planar domains  with $\partial
\Sigma_i\subset
\partial B_{R_i}$ where $R_i\to \infty$, then there is a subsequence with either
$\cSu = \emptyset$ or $\cSt = \emptyset$.
\end{theorem}

The case where $\cSt = \emptyset$ was dealt with in the previous
section.

Common for both the ULSC case and the case where $\cSu$ is empty
is that the limits are always  laminations by flat parallel planes
and  the singular sets are always  closed subsets contained in the
union of the planes.  This is the content of the next theorem:

\begin{theorem} \label{t:tab}
\cite{CM7}
 Let $\Sigma_i \subset B_{R_i}=B_{R_i}(0)\subset \RR^3$
be a sequence of compact embedded minimal surfaces with fixed
genus and with $\partial \Sigma_i\subset \partial B_{R_i}$ where
$R_i\to \infty$.   If
\begin{equation}
\sup_{B_1\cap \Sigma_i}|A|^2\to \infty \, ,
 \end{equation}
   then there exists a subsequence $\Sigma_j$,
     a lamination $\cL=\{x_3=t\}_{ \{ t \in \cI \} }$ of $\RR^3$
 by parallel
 planes (where $\cI \subset \RR$ is a closed set), and a closed nonempty set
 $\cS$ in the union of the leaves of $\cL$ such that
 after a rotation of $\RR^3$:
 \begin{enumerate}
\item[(A)] For each $1>\alpha>0$, $\Sigma_j\setminus \cS$
converges in the $C^{\alpha}$-topology to the lamination $\cL
\setminus \cS$.
\item[(B)]  $\sup_{B_{r}(x)\cap \Sigma_j}|A|^2\to\infty$ as $j \to \infty$ for
all $r>0$ and $x\in \cS$.  (The curvatures blow up along $\cS$.)
\end{enumerate}
\end{theorem}

One can get both types of curvature blow-up by considering the
family of embedded minimal planar domains known as the Riemann
examples.
 Modulo translations and rotations, this is a
two-parameter family of periodic minimal surfaces, where the
parameters can be thought of as the size of the necks and the
angle  from one fundamental domain to the next. By choosing the
two parameters appropriately, one can produce sequences of Riemann
examples that illustrate both of the two structure theorems (cf.
figures \ref{f:01a} and \ref{f:01b}):
\begin{enumerate}
\item
If we take a sequence of Riemann examples where the neck size is
fixed and the angles go to $\frac{\pi}{2}$, then the surfaces with
angle near $\frac{\pi}{2}$ can be obtained by gluing together two
oppositely--oriented double spiral staircases. Each double spiral
staircase looks like a helicoid. This sequence of Riemann examples
converges to a foliation by parallel planes.  The convergence is
smooth away from the axes of the two helicoids (these two axes are
the singular set $\cS$ where the curvature blows up). The sequence
is ULSC since the size of the necks is fixed and  thus illustrates
the first structure theorem, Theorem \ref{t:t5.1}.
\item
If we take a sequence of examples where the neck sizes go to zero,
then we get a sequence  that is {\it{not}} ULSC.  However, the
surfaces can be cut along   short curves into collections of
graphical pairs of pants.  The short curves converge to points and
the graphical pieces converge to flat planes except at these
points, illustrating the second structure theorem, Theorem
\ref{t:t5.2a} below.
\end{enumerate}

With these examples in mind, we are now ready to state our second
main structure theorem describing the case where $\cSu$ is empty.

\begin{theorem} \label{t:t5.2a}
 Let  a sequence $\Sigma_i$, limit lamination $\cL$, and singular
 set $\cS$
be  as in Theorem \ref{t:tab}.  If $\cSu = \emptyset$ and
\begin{equation}
\sup_{B_1\cap \Sigma_i}|A|^2\to \infty \, ,
 \end{equation}
   then $\cS=\cSt$ by \eqr{e:} and
 \begin{enumerate}
\item[($C_{neck}$)]
Each point $y$ in $\cS$ comes with a sequence of
{\underline{graphs}} in $\Sigma_j$ that converge  to the plane $\{
x_3 = x_3 (y) \}$. The convergence is in the $C^{\infty}$ topology
away from the point $y$ and possibly also one other point  in $\{
x_3 = x_3 (y) \} \cap \cS$.  If the convergence is   away from one
point, then these graphs are defined over annuli; if the
convergence is away from two points, then the graphs are defined
over disks with two subdisks removed.
\end{enumerate}
\end{theorem}

\part{Global theory of minimal surfaces in $\RR^3$}

\stepcounter{section}

\subsection{Minimal surfaces with finite total curvature}

A complete, non-compact, minimal immersion of a surface $\Sigma$
in $\RR^3$ is said to have {\it finite total curvature}  (or ftc)
if
\begin{equation}    \label{e:ftc}
    \int_{\Sigma} |A|^2 < \infty \, .
\end{equation}
The simplest examples are the plane (where the total  curvature is
zero) and the catenoid.  In the case of the catenoid,   the Gauss
map gives a conformal diffeomorphism to the sphere punctured at
the north and south poles.  Since  $A$ is the differential of the
Gauss map, it follows that the catenoid has total curvature
$8\pi$.

The fundamental result for minimal surfaces with ftc is that they
are all {\underline{conformally}} diffeomorphic
  to compact Riemann surfaces with a finite number of
points removed.  Namely, we have the following result of R.
Osserman (see \cite{Os}):

\begin{theorem}     \label{t:oss}
Let $\Sigma \subset \RR^3$ be a complete minimal immersion with
finite total curvature.  Then:
\begin{enumerate}
\item  $\Sigma$ is conformally diffeomorphic to a compact Riemann
surface with a finite set of points removed.  Each point
corresponds to an end of the surface.
\item  $\Sigma$ is proper.
\item  The Weierstrass data extends across the punctures
meromorphically.  (See Section \ref{s:21} for the definition of
the Weierstrass representation.)
\item  The total curvature is an integral multiple of $8\pi$.
\end{enumerate}
\end{theorem}

One application of this theorem is a classification of the ends of
an embedded minimal surface with finite total curvature.  Namely,
one can show that any such end is asymptotic to either a plane or
to half of a catenoid.

\subsection{The uniqueness of the catenoid}

Given the structure result, Theorem \ref{t:oss}, it is natural to
try to understand the space of  minimal surfaces with finite total
curvature in terms of the genus and the number of ends.  Two such
results were proven by R. Schoen and F. Lopez--A. Ros,
respectively. The theorem of Schoen says that the catenoid is the
unique embedded minimal surface with finite total curvature and
exactly two ends.

\begin{theorem}    \cite{Sc2}   \label{t:sc2}
If $\Sigma \subset \RR^3$ is a connected embedded minimal surface
with finite total curvature and exactly two ends, then $\Sigma$ is
a catenoid.
\end{theorem}

The theorem of Lopez and Ros  says that the catenoid and the plane
are the only embedded minimal surfaces with finite total curvature
and genus zero.  Here, ``genus zero'' means conformal to the
sphere with a finite set of points removed.

\begin{theorem}    \cite{LoRo}   \label{t:loro}
If $\Sigma \subset \RR^3$ is an embedded minimal surface with
finite total curvature and genus zero, then $\Sigma$ is a catenoid
or a plane.
\end{theorem}

\section{Global theory of embedded minimal surfaces}      \label{s:20}

Recent years have seen breakthroughs on many long--standing
problems in the global theory of minimal surfaces in $\RR^3$. This
is an enormous subject and, rather than give a comprehensive
treatment, we will mention a few important results which fit well
with the rest of this survey. Throughout this section, $\Sigma$
will be a complete properly embedded minimal surface in $\RR^3$
(recall that properness here means that the intersection of
$\Sigma$ with any compact subset of $\RR^3$ is compact).

We say that $\Sigma$ has finite topology if it is homeomorphic to
a closed Riemann surface with a finite number of punctures; the
genus of $\Sigma$ is then the genus of this Riemann surface and
the number of punctures is the number of ends.  It follows that a
neighborhood of each puncture corresponds to a properly embedded
annular end of $\Sigma$.  Perhaps surprisingly at first, the more
restrictive case is when $\Sigma$ has more than one end. The
reason for this is that a barrier argument gives a stable minimal
surface between any pair of ends.  Such a stable surface is then
asymptotic to a plane (or catenoid), essentially forcing each end
to live in a half--space. Using this restriction, P. Collin
proved:

\begin{theorem}  \cite{Co}      \label{t:co}
Each end of a complete properly embedded minimal surface with
finite topology and at least two ends is asymptotic to a plane or
catenoid.
\end{theorem}

In particular, outside some compact set, $\Sigma$ is given by a
finite collection of disjoint graphs over a common plane (and has
finite total curvature).

As mentioned above, Collin essentially proved Theorem \ref{t:co}
by showing that an embedded annular end that lives in a
half--space must have finite total curvature.  \cite{CM14} used
the one--sided curvature estimate to strengthen this from a
half--space to a strictly larger  cone, and in the process give a
very different proof of Collin's theorem.

\begin{theorem}  \cite{CM14}      \label{t:coCM}
There exists $\epsilon > 0$ so that any complete properly embedded
minimal annular end contained in the cone
\begin{equation}
    \{ x_3 \geq - \epsilon \, (x_1^2 + x_2^2 + x_3^2)^{1/2}  \}
\end{equation}
 is asymptotic to a plane or catenoid.
\end{theorem}

\vskip2mm When $\Sigma$ has only one end (e.g., for the helicoid),
it need not have finite total curvature so the situation is more
delicate.  However, the regularity results of the previous section
can be applied.  For example, if $\Sigma$ is a (non--planar)
embedded minimal disk, then we get a multi--valued graph structure
 away from a ``one--dimensional singular set.''  Using
Theorems \ref{t:t0.1} and \ref{t:t2}, W. Meeks and H. Rosenberg
proved the uniqueness of the helicoid:

\begin{theorem}     \cite{MeRo1}     \label{t:mero}
The plane and helicoid are the only complete properly embedded
simply--connected minimal surfaces in $\RR^3$.
\end{theorem}

This uniqueness has many applications. Recall that if we take a
sequence of rescalings of the helicoid, then the singular set
$\cS$ for the convergence is the vertical axis perpendicular to
the leaves of the foliation. In \cite{Me1}, W. Meeks used this
fact together with the uniqueness of the helicoid to prove that
the singular set $\cS$ in Theorem \ref{t:t0.1} is always a
straight line perpendicular to the foliation.  Recently, W. Meeks
and M. Weber have constructed a {\it local} example (i.e., a
sequence of embedded minimal disks in a unit ball) where $\cS$ is
a circle.

\vskip2mm We have not even touched on the case where $\Sigma$ has
infinite topology (e.g., when $\Sigma$ is one of the Riemann
examples). This is an area of much current research, see
\cite{CM5},  the work of Meeks, J. Perez and A. Ros,
\cite{MePRs1}, \cite{MePRs2}, and \cite{MePRs3}, the survey
\cite{MeP} and references therein.

\vskip2mm We close this section with a local analog of the
two--ended case.  Namely, in \cite{CM9}, we proved
 that any embedded minimal annulus in a ball
(with boundary in the boundary of the ball and) with a small neck
can be decomposed by a simple closed geodesic into two graphical
sub--annuli.  Moreover, we gave a sharp bound for the length of
this closed geodesic in terms of the separation (or height)
between the graphical sub--annuli.  This  serves to illustrate our
``pair of pants'' decomposition from \cite{CM5} in the special
case where the embedded minimal planar domain is an annulus  (we
will not touch on this further here). The catenoid
\begin{equation}
    \{x_1^2+x_2^2=\cosh^2 x_3\}
\end{equation}
 is the prime example of an embedded minimal annulus.

The precise statement of this decomposition for annuli is:

\begin{theorem} \cite{CM9} \label{t:nitsche}
There exist $\epsilon>0$, $C_1 ,\, C_2,\, C_3>1$ so:   If
$\Sigma\subset B_{R}\subset \RR^3$ is an embedded minimal annulus
with $\partial \Sigma\subset \partial B_{R}$ and
$\pi_1(B_{\epsilon R}\cap \Sigma)\ne 0$, then there exists a
simple closed geodesic $\gamma \subset \Sigma$ of length $\ell$ so
that:

\begin{itemize}
\item
 The curve $\gamma$ splits the connected component of $B_{R/C_1}\cap
\Sigma$ containing it into
 two annuli
$\Sigma^{+} , \Sigma^{-}$ each with $\int |A|^2 \leq 5 \, \pi$.
 \item
 Furthermore, $\Sigma^{\pm}
\setminus \cT_{C_2 \, \ell}(\gamma)$ are graphs with gradient
$\leq 1$.
\item
Finally, $\ell \log (R/\ell) \leq C_3\,h$ where the separation $h$
is given by \begin{equation}
 h=\min_{x_{\pm} \in\partial B_{
R/C_1 }\cap \Sigma^{\pm}}|x_+ -x_- | \, .
\end{equation}
\end{itemize}
\end{theorem}

Here  $\cT_{s}(S) \subset \Sigma$ denotes the intrinsic
$s$--tubular neighborhood of a subset $S\subset \Sigma$.

\section{The Calabi-Yau conjectures}

Recall that an immersed submanifold in $\RR^n$   is {\it proper}
if the pre-image of any compact subset of $\RR^n$ is compact in
the surface.  This property has played an important role in the
theory of minimal submanifolds and  many of the classical theorems
in the subject assume that the submanifold is proper.

It is easy to see that any compact submanifold is automatically
proper.  On the other hand, there is no reason to expect a general
immersion (or even embedding) to be proper. For example, the
non-compact curve parametrized in polar coordinates by
\begin{equation}
    \rho (t ) = \pi + \arctan (t) \, , \, \theta (t) = t
\end{equation}
spirals infinitely between the circles of radius $\pi/2$ and
$3\pi/2$. However, it was long thought that a minimal immersion
(or embedding) should be better behaved.  This principle was
captured by the Calabi-Yau conjectures, dating back to the 1960s.
Much work has been done on them over the past four decades. Their
original form  was given in 1965 in \cite{Ca} where E. Calabi made
the following two conjectures about minimal surfaces (see also
S.S. Chern, page 212 of \cite{Ch} and S.T. Yau's 1982 problem
list):

\begin{conjecture} \label{con:1}
``Prove that a complete  minimal hypersurface in $\RR^n$ must be
unbounded.''
\end{conjecture}

Calabi continued: ``It is known that there are no compact minimal
submanifolds of $\RR^n$ (or of any simply connected complete
Riemannian manifold with sectional curvature $\leq 0$). A more
ambitious conjecture is'':

\begin{conjecture} \label{con:2}
``A complete [nonflat] minimal hypersurface in $\RR^n$ has an
unbounded projection in every $(n-2)$--dimensional flat
subspace.''
\end{conjecture}

 The \underline{immersed}  versions of these conjectures turned out to be
 false.  As mentioned above,  Jorge and Xavier, \cite{JXa2},  constructed
nonflat minimal immersions contained between two parallel planes
in 1980, giving a counter-example to the immersed version of the
more ambitious Conjecture \ref{con:2}.
 Another significant development came in 1996,  when
 N. Nadirashvili, \cite{Na1}, constructed a complete
immersion of a minimal disk into the unit ball in $\RR^3$, showing
that Conjecture \ref{con:1} also failed for immersed surfaces; see
  \cite{MaMo1}, \cite{LMaMo1}, \cite{LMaMo2},
 for other topological types than disks.

 It is clear from the definition of proper that a proper
minimal surface in $\RR^3$ must be unbounded, so the examples of
Nadirashvili are not proper.  Much less obvious is that the plane
is the only complete {\underline{proper}} immersed minimal surface
in a halfspace.  This is however a consequence of the strong
halfspace theorem of D. Hoffman and W. Meeks:

\begin{theorem}
\cite{HoMe} A complete connected {\underline{properly}} immersed
minimal surface $\Sigma \subset \{ x_3 > 0 \} \subset \RR^3$ must
be a horizontal plane $\{ x_3 = {\text{Constant}} \}$.
\end{theorem}

The main result of \cite{CM16} is an effective version of
properness for disks, giving a chord arc bound.      Obviously,
intrinsic distances are larger than extrinsic distances, so the
significance of a chord arc bound is the reverse inequality, i.e.,
a bound on intrinsic distances from above by extrinsic distances.
This is accomplished in the next theorem:

\begin{theorem}     \label{t:1}
\cite{CM16} There exists a  constant  $C > 0$  so that if $\Sigma
\subset \RR^3$ is an embedded minimal disk, $\cB_{2R}=\cB_{2R}(0)$
is an intrinsic ball in $\Sigma \setminus
\partial \Sigma$ of radius $2R$, and
$ \sup_{\cB_{r_0}}|A|^2>r_0^{-2}$ where $R>  r_0$,
 then for $x \in \cB_R$
 \begin{equation}   \label{e:t1}
    C \, \dist_{\Sigma}(x,0)<   |x| + r_0   \, .
 \end{equation}
\end{theorem}

The assumption of a lower curvature bound, $
\sup_{\cB_{r_0}}|A|^2>r_0^{-2}$, in the theorem is a necessary
normalization for a chord arc bound.  This can easily be seen by
rescaling and translating the helicoid.

 Properness of a complete embedded minimal disk   is an immediate
consequence of Theorem \ref{t:1}.  Namely, by \eqr{e:t1}, as
intrinsic distances go to infinity, so do extrinsic distances.
Precisely, if $\Sigma$ is flat, and hence a plane, then obviously
$\Sigma$ is proper and if it is non-flat, then
$\sup_{\cB_{r_0}}|A|^2>r_0^{-2}$ for some $r_0>0$ and hence
$\Sigma$ is proper by \eqr{e:t1}.

A consequence of Theorem \ref{t:1} together with the one-sided
curvature estimate of \cite{CM6} (i.e., theorem $0.2$ in
\cite{CM6}) is the following  version of that estimate for
intrinsic balls:

\begin{corollary}  \label{t:one-sided}
\cite{CM16} There exists $\epsilon>0$, so that if
\begin{equation}
    \Sigma \subset \{x_3>0\} \subset \RR^3
\end{equation}
 is an embedded
minimal disk with $\cB_{2R} (x) \subset \Sigma \setminus
\partial \Sigma$ and $|x|<\epsilon\,R$,  then
\begin{equation}        \label{e:graph2}
\sup_{ \cB_{R}(x) } |A_{\Sigma}|^2 \leq R^{-2} \, .
\end{equation}
\end{corollary}

 As a corollary of this intrinsic one-sided curvature estimate we get that
the second, and more ambitious, of Calabi's conjectures is also
true for {\underline{embedded}} minimal disks.

\vskip2mm   In fact,  \cite{CM16}  proved both of Calabi's
conjectures and properness also for embedded surfaces with finite
topology. Recall that a surface $\Sigma$ is said to have finite
topology if it is homeomorphic to a closed Riemann surface with a
finite set of points removed or ``punctures''. Each puncture
corresponds to an end of $\Sigma$.

The following generalization of the halfspace theorem gives
Calabi's second, more ambitious, conjecture for embedded surfaces
with finite topology:

\begin{theorem}     \label{c:finite}
The plane is the only complete embedded minimal surface with
finite topology in $\RR^3$
 in a halfspace.
\end{theorem}

Likewise, we get the properness of embedded surfaces with finite
topology:

\begin{theorem}     \label{c:3rd}
A complete embedded minimal surface with finite topology in
$\RR^3$ must be proper.
\end{theorem}

Compare also W. Meeks and H. Rosenberg, \cite{MeRo2}.

\vskip2mm There has been extensive  work on both properness   and
the halfspace property     assuming various
{\underline{curvature}} {\underline{bounds}}. Jorge and Xavier,
\cite{JXa1} and \cite{JXa2}, showed that there cannot exist a
complete immersed minimal surface with {\underline{bounded}}
{\underline{curvature}} in $\cap_i \{ x_i > 0 \}$; later Xavier
proved that the plane is the only such surface in a halfspace,
\cite{Xa}. Recently, G.P. Bessa, Jorge and G. Oliveira-Filho,
\cite{BJO}, and H. Rosenberg, \cite{Ro3}, have shown that if such
a surface is embedded,  then it must be proper. This properness
was extended to embedded minimal surfaces with locally bounded
curvature and finite topology by Meeks and Rosenberg in
\cite{MeRo1}; finite topology was subsequently replaced by finite
genus in \cite{MePRs1} by Meeks, J. Perez and A. Ros.

Inspired by Nadirashvili's examples, F. Martin and S. Morales
 constructed  in \cite{MaMo2} a complete bounded minimal
immersion which is proper in the (open) unit ball.  That is, the
preimages of compact subsets of the (open) unit ball are compact
in the surface and the image of the surface accumulates on the
boundary of the unit ball. They extended this  in
 \cite{MaMo3} to show that any convex, possibly noncompact or
nonsmooth, region of $\RR^3$ admits a proper complete minimal
immersion of the unit disk; cf.  \cite{Na2}.

\part{Constructing minimal surfaces}   \label{p:4}

Thus far, we have mainly dealt with regularity and {\it a priori}
estimates but have ignored questions of existence.  In this part
we survey some of the most useful existence results for minimal
surfaces. Section \ref{s:plateau}   gives an overview of the
classical Plateau problem.  Section \ref{s:21} recalls the
classical Weierstrass representation, including a few modern
applications, and the Kapouleas desingularization method.  Section
\ref{s:22} deals with producing area minimizing surfaces  and
questions of embeddedness. Finally, Section \ref{s:24} recalls the
min--max construction for producing unstable minimal surfaces and,
in particular, doing so while controlling the topology and
guaranteeing embeddedness.

\section{The Plateau Problem}  \label{s:plateau}

The following fundamental existence problem for minimal surfaces
is known as the {\it Plateau problem}:   Given a closed curve
$\Gamma$, find a minimal surface with boundary $\Gamma$. There are
various solutions to this problem depending on the exact
definition of a surface (parameterized disk, integral current,
$\ZZ_2$ current, or rectifiable varifold).
 We shall consider the
version of the Plateau problem for parameterized disks; this was
solved independently by J. Douglas   and T. Rado. The
generalization to Riemannian manifolds is due to C. B. Morrey.

\begin{theorem}     \label{t:plat}
Let $\Gamma \subset \RR^3$ be a piecewise $C^1$ closed Jordan
curve. Then there exists a piecewise $C^1$ map $u$ from $D \subset
\RR^2$ to $\RR^3$ with $u(\partial D) \subset \Gamma$ such that
the image  minimizes area among all disks with boundary $\Gamma$.
\end{theorem}

The solution $u$ to the Plateau problem  above can easily be seen
to be a branched conformal immersion.
 R. Osserman proved that $u$ does not have true interior branch
 points; subsequently,
R. Gulliver and W. Alt showed that $u$ cannot have false branch
points either.

Furthermore, the solution $u$ is as smooth as the boundary curve,
even up to the boundary.    A very general version of this
boundary regularity was proven by S. Hildebrandt; for the case of
surfaces in $\RR^3$, recall the following result of J. C. C.
Nitsche:

\begin{theorem}              \label{t:nitsche2}
If $\Gamma$ is a regular Jordan curve of class $C^{k,\alpha}$
where $k \geq 1$ and $0 < \alpha < 1$, then a solution $u$ of the
Plateau problem is $C^{k,\alpha}$ on all of $\bar{D}$.
\end{theorem}

The optimal boundary regularity theorem in higher dimensions was
proven by R. Hardt and L. Simon in \cite{HaSi}.

\section{The Weierstrass representation}        \label{s:21}

The classical {\it Weierstrass representation} (see \cite{HoK} or
\cite{Os}) takes holomorphic data (a Riemann surface, a
meromorphic function, and a holomorphic one--form) and associates
to these data a minimal surface in $\RR^3$.
 To be precise, given
a Riemann surface $\Omega$, a meromorphic function $g$ on
$\Omega$, and a holomorphic one--form $\phi$ on $\Omega$, then we
get  a (branched) conformal minimal immersion $F: \Omega \to
\RR^3$ by
\begin{equation}    \label{e:ws1}
    F(z) = {\text{Re }} \int_{\zeta \in \gamma_{z_0,z}}
\left( \frac{1}{2} \, (g^{-1} (\zeta) - g (\zeta) )
    , \frac{i}{2} \, (g^{-1} (\zeta)
    +g (\zeta) ) , 1 \right) \, \phi (\zeta) \, .
\end{equation}
Here $z_0 \in \Omega$ is a fixed base point and the integration is
along a path $\gamma_{z_0,z}$ from $z_0$ to $z$. The choice of
$z_0$ changes $F$ by adding a constant. In general, the map $F$
may depend on the choice of path (and hence may not be
well--defined);  this is known as ``the period problem'' (see M.
Weber and M. Wolf, \cite{WeWo}, for the latest developments).
However, when $g$ has no zeros or poles and $\Omega$ is simply
connected, then
 $F(z)$ does not depend on the choice of path
$\gamma_{z_0,z}$.

Three standard constructions of minimal surfaces from Weierstrass
data are
 \begin{align}
 &g (z) = z, \, \phi (z)  = dz/z ,
\, \Omega = \CC \setminus \{ 0 \} {\text{ giving a catenoid}}
\, , \\
&g (z) = z, \, \phi (z)  = z\, dz , \, \Omega = \CC   {\text{
giving Enneper's surface}}
\, , \label{e:enneper} \\
  &g (z) = \e^{iz} , \,
\phi (z) = dz , \, \Omega = \CC  {\text{ giving a helicoid}} \, .
\label{e:hel}
\end{align}

The Weierstrass representation is particularly useful for
constructing immersed minimal surfaces.  For example, in
\cite{Na1},  N. Nadirashvili used it to construct a complete
immersed minimal surface in the unit ball in $\RR^3$ (see also
\cite{JXa2} for the case of a slab). In particular, Nadirashvili's
surface is not proper, i.e., the intersections with compact sets
are not necessarily compact.

 Typically, it is rather difficult to prove that the resulting
immersion is an embedding (i.e., is $1$--$1$), although there are
some interesting cases where this can be done.    The first modern
example was \cite{HoMe} where D. Hoffman and W. Meeks proved that
the surface constructed by Costa was embedded; this was the first
new complete finite topology properly embedded minimal surface
discovered since the classical catenoid, helicoid, and plane. This
led to the discovery of many more such surfaces
 (see \cite{HoK} and \cite{Ro1} for more discussion).

Very recently in \cite{HoWeWo1} (see also \cite{HoWeWo2}), D.
Hoffman, M. Weber, and M. Wolf have used the Weierstrass
representation to construct a genus one properly embedded minimal
surface asymptotic to the helicoid. They construct this ``genus
one helicoid''   as the limit of a continuous one-parameter family
of screw-motion invariant minimal surfaces--also asymptotic to the
helicoid--that have genus equal to one in the quotient.

In \cite{CM11}, we used the Weierstrass representation to
construct a sequence of embedded minimal disks
\begin{equation}
    \Sigma_i \subset B_1 = B_1 (0)\subset \RR^3
\end{equation}
 with $\partial \Sigma_i \subset
\partial B_1$ where
 the curvatures blow up only at $0$  and $\Sigma_i \setminus \{ {\text{$x_3$--axis}} \}$
consists of two multi--valued graphs for each $i$. Furthermore,
$\Sigma_i  \setminus \{ x_3 = 0 \}$ converges to two
  embedded minimal disks $\Sigma^- \subset \{ x_3 < 0 \}$
and  $\Sigma^+ \subset \{ x_3 > 0 \}$ each of which spirals
 into $\{ x_3 =
0 \}$ and thus is not proper.  (This should be contrasted with
Theorem \ref{t:t0.1} where the {\it complete} limits are planes
and hence proper.)

\vskip2mm N. Kapouleas has developed another method to construct
  complete embedded minimal surfaces with finite total
curvature.  For instance, in \cite{Ka3}, he shows that (most)
collections of coaxial   catenoids and planes can be
desingularized to get complete embedded minimal surfaces with
finite total curvature. The Costa surface above had genus one and
three ends (that is to
 say, it is homeomorphic to a torus with three punctures).  In the
 Kapouleas construction, one could start with a plane and catenoid
 intersecting in a circle and then desingularize this circle using
 suitably scaled and bent Scherk surfaces
 to get a finite genus embedded surface with three
 ends.  (This desingularization process adds handles, i.e.,
 increases the genus.)  In this manner, Kapouleas gets an enormous
 number of new examples; see also the gluing construction of S.D. Yang,
 \cite{Y}, which uses catenoid necks to glue together nearby
 minimal surfaces.

\section{Area--minimizing surfaces}     \label{s:22}

Perhaps the most natural way to construct minimal surfaces is to
look for ones which minimize area, e.g., with fixed boundary, or
in a homotopy class, etc.  This has the advantage that often it is
possible to show that the resulting surface is embedded.  We
mention a few results along these lines.

The first embeddedness result, due to   W. Meeks and S.T. Yau,
shows that if the boundary curve is embedded and lies on the
boundary of a smooth mean convex set (and it is null--homotopic in
this  set), then it bounds an embedded least area disk.

\begin{theorem}  \cite{MeYa1}             \label{t:my1}
Let $M^3$ be a compact Riemannian three--manifold whose boundary
is mean convex and let $\gamma$ be a simple closed curve in
$\partial M$ which is null--homotopic in $M$; then $\gamma$ is
bounded by a least area disk and any such least area disk is
properly embedded.
\end{theorem}

Note that some restriction on the boundary curve $\gamma$ is
certainly necessary. For instance, if the boundary curve was
knotted (e.g., the trefoil), then it could not be spanned by any
embedded disk (minimal or otherwise).  Prior to the work of Meeks
and Yau, embeddedness was known for extremal boundary curves in
$\RR^3$ with small total curvature by the work of R. Gulliver and
J. Spruck \cite{GuSp}; see chapter $4$ in \cite{CM1} for other
results and further discussion.  Recently, in \cite{EWWi}, T.
Ekholm, B. White, and D. Wienholtz proved that minimal surfaces
whose boundary has total curvature less than $4\pi$ also must be
embedded.

 If we instead fix a homotopy class of maps, then the two fundamental existence results
 are due to J. Sacks--K. Uhlenbeck and R. Schoen-- S.T. Yau (with embeddedness proven
 by W. Meeks--S.T. Yau and M. Freedman--J. Hass--P. Scott,
 respectively):

\begin{theorem}     \label{p:existence1} \cite{SaUh}, \cite{MeYa2}
Given $M^3$, there exist conformal  (stable) minimal immersions
$u_1 , \dots , u_m  : \SS^2 \to M$ which generate $\pi_2 (M)$ as a
$\ZZ[\pi_1 (M)]$ module.  Furthermore,
\begin{itemize}
\item
If $u: \SS^2 \to M$ and $[u]_{\pi_2} \ne 0$, then $\Area (u) \geq
\min_i \Area (u_i)$.
 \item Each $u_i$ is either an
embedding or a $2$--$1$ map onto an embedded $2$--sided $\RP^2$.
\end{itemize}
\end{theorem}

\begin{theorem}     \label{p:existence2} \cite{ScYa2}, \cite{FHS}
If $\Sigma^2$ is a closed surface with genus $g>0$ and $i_0 :
\Sigma \to M^3$ is an embedding which induces an injective map on
$\pi_1$, then there is a least area embedding with the same action
on $\pi_1$.
\end{theorem}

In \cite{MeSiYa}, W. Meeks, L. Simon, and S.T. Yau find an
embedded sphere minimizing area in an isotopy class in a closed
$3$--manifold.

\vskip2mm We end this section by mentioning two applications of
Theorem \ref{p:existence2}.  First, in \cite{CM13}, we showed that
any topological $3$--manifold $M$ had an open set of metrics so
that, for each such metric, there was a sequence of embedded
minimal tori whose area went to infinity.  In \cite{De}, B. Dean
showed that this was true for every genus $g \geq 1$.

\section{The min--max construction of minimal surfaces} \label{s:24}

Variational arguments can also be used to construct higher index
(i.e., non--minimizing) minimal surfaces using   the topology of
the space of surfaces. There are two basic approaches:
\begin{itemize}
\item
Applying Morse theory to the energy functional on the space of
maps from a fixed surface $\Sigma$ to $M$.
\item
Doing a min--max argument over  families of (topologically
non--trivial) sweep--outs of $M$.
\end{itemize}
The first approach has the advantage that the topological type of
the minimal surface is easily fixed; however, the second approach
has been  more successful at producing embedded minimal surfaces.
We will highlight a few key  results below but refer to \cite{CD}
for a thorough treatment.

Unfortunately, one cannot directly apply Morse theory to the
energy functional on the space of maps from a fixed surface
because of a lack of compactness (the Palais--Smale Condition C
does not hold).  To get around this difficulty, J. Sacks and K
Uhlenbeck, \cite{SaUh}, introduce a family of perturbed energy
functionals which do satisfy Condition C and then obtain   minimal
surfaces as limits of critical points for the perturbed problems:

\begin{theorem} \cite{SaUh}     \label{t:nonasph}
If   $\pi_k (M) \ne 0$ for some $k>1$, then there exists a
branched immersed minimal $2$--sphere in $M$ (for any metric).
\end{theorem}

This was sharpened somewhat by M. Micallef and D. Moore,
\cite{MiMo}, (showing that the index of the minimal sphere was at
most $k-2$), who used it to prove a generalization of the sphere
theorem.  See A. Fraser, \cite{Fr}, for a generalization to a free
boundary problem.

\begin{figure}[htbp]
    \setlength{\captionindent}{4pt}
    %\begin{minipage}[t]{0.5\textwidth}
    \centering\input{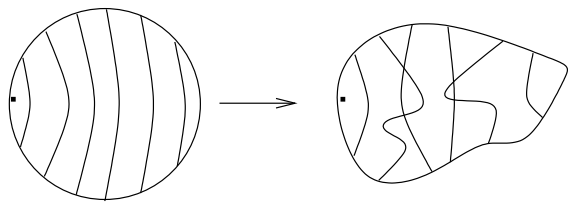}
    \caption{A $1$--parameter family of curves on a $2$--sphere which
induces a map $F:\SS^2 \to
    \SS^2$ of degree $1$. First published in Surveys in Differential Geometry,
    volume IX, in 2004, published by International Press.}\label{f:fsweep}
    %\end{minipage}
\end{figure}

\vskip2mm The basic idea of constructing minimal surfaces via
min--max arguments and sweep--outs goes back to G. Birkhoff, who
developed it to construct simple closed geodesics on spheres.  In
particular, when $M$ is a topological $2$--sphere, we can find a
$1$--parameter family of curves starting and ending at point
curves so that the induced map $F:\SS^2 \to \SS^2$ (see figure
\ref{f:fsweep}) has nonzero degree.    The min--max argument
produces a nontrivial closed geodesic of length less than or equal
to the longest curve in the initial one--parameter family.  A
curve shortening argument gives that the geodesic obtained in this
way is simple.

In \cite{Pi}, J. Pitts applied a similar argument and   geometric
measure theory to get that every closed Riemannian three manifold
has an embedded minimal surface (his argument was for dimensions
up to seven), but he did not estimate the genus of the resulting
surface.  Finally, F. Smith (under the direction of L. Simon)
proved (see \cite{CD}):

\begin{theorem} \cite{Sm}   \label{t:smith}
Every metric   on  a topological $3$--sphere $M$ admits an
embedded minimal $2$--sphere.
\end{theorem}

The main new contribution of Smith was to control the topological
type of the resulting minimal surface while keeping it embedded;
see also J. Pitts and J.H. Rubinstein, \cite{PiRu}, for some
generalizations.

\part{Some applications of minimal surfaces}

\stepcounter{section}

In this last part, we discuss very briefly a few applications of
minimal surfaces.  As mentioned in the introduction, there are
many to choose from and we have selected just a few. See
\cite{Am3}, \cite{Fi}, \cite{HaKLi}, \cite{Ho}, \cite{Ta}, and
\cite{Th} for other applications.

\subsection{The positive mass theorem}  \label{s:posmass}

The (Riemannian version of the) {\it positive mass theorem} of R.
Schoen and S.T. Yau states that an asymptotically flat
$3$-manifold $M$ with non-negative scalar curvature must have
positive mass. The Riemannian manifold $M$ here arises as a
maximal space-like slice in a $3+1$-dimensional space-time
solution of Einstein's equations.

The asymptotic flatness of $M$ comes from that the space-time
models an isolated gravitational system and hence is a
perturbation of the vacuum solution outside a large compact set.
 To make this precise, suppose for simplicity that $M$ has only
 one end;
 $M$ is then said to be {\it asymptotically flat} if there is a compact
set $\Omega \subset M$ so that $M \setminus \Omega$ is
diffeomorphic to $\RR^3 \setminus B_R (0)$ and the metric on $M
\setminus \Omega$ can be written as
\begin{equation}    \label{e:pdef}
    g_{ij} = \left(1 + \frac{\cM}{ 2 \, |x|} \right)^4 \, \delta_{ij} + p_{ij} \, ,
\end{equation}
where
\begin{equation}    \label{e:pbds}
    |x|^2 \, |p_{ij}| + |x|^3 \, |D \, p_{ij}|
    + |x|^4 \, |D^2 \, p_{ij}|  \leq C \, .
\end{equation}
The constant $\cM$ is the so called {\it mass} of $M$.  Observe
that the metric $g_{ij}$ is a perturbation of the metric on a
constant-time slice in the Schwarzschild space-time of mass $\cM$;
that is to say, the Schwarzschild metric has $p_{ij} \equiv 0$.

A tensor $h$ is said to be $O(|x|^{-p})$ if $|x|^p \, |h| +
|x|^{p+1} \, |D \, h| \leq C$.  For example, an easy calculation
shows that
\begin{align}   \label{e:pdefst}
    g_{ij} &= \left(1 + 2 \, \cM / |x| \right) \, \delta_{ij} +O(|x|^{-2}) \, , \\
    \sqrt{g} &\equiv \sqrt{\det g_{ij} } = 1 + 3 \, \cM \, |x|^{-1} +
    O(|x|^{-2}) \, . \notag
\end{align}

The positive mass theorem states that the mass $\cM$ of such an
$M$ must be non-negative:

\begin{theorem} \label{t:pm1}
\cite{ScYa1}  With $M$ as above, $\cM \geq 0$.
\end{theorem}

There is a rigidity theorem as well which states that the mass
vanishes only when $M$ is isometric to $\RR^3$:

\begin{theorem} \label{t:pm2}
\cite{ScYa1} If $|\nabla^3 p_{ij}|  = O(|x|^{-5})$
 and  $\cM = 0$ in Theorem \ref{t:pm1},
then $M$ is {\underline{isometric}} to $\RR^3$.
\end{theorem}

 We will give a very brief overview of the proof of
 Theorem \ref{t:pm1}, showing in the process where minimal surfaces appear.

\begin{proof}
(Sketch.) The argument will by contradiction, so suppose that the
mass is negative.  Note first that   a ``rounding off'' argument
shows that the metric on $M$ can be perturbed to have
{\underline{positive}} scalar curvature outside of a compact set
and still have negative mass.

It is not hard to prove that the slab between two parallel planes
is mean-convex.  That is, we have the following:

\begin{lemma} \label{l:ecalc}
If $\cM < 0$ and $M$ is asymptotically flat, then there exist $R_0
, h > 0$ so that for $r > R_0$ the sets
\begin{equation}
    C_r = \{ |x|^2  \leq r^2  \, , \, -h \leq x_3 \leq h \}
\end{equation}
have strictly mean convex boundary.
\end{lemma}

Since the compact set $C_r$ is mean convex, we can solve the
Plateau problem (as in Section \ref{s:plateau}) to get an area
minimizing (and hence stable) surface $\Gamma_r \subset C_r$ with
boundary
\begin{equation}
\partial \Gamma_r = \{
|x|^2 = r^2 \, , \, x_3 = h \} \, .
\end{equation}
  Using the disk
$\{ |x|^2 \leq r^2 \, , x_3 = h \}$ as a comparison surface, we
get uniform local area bounds for any such $\Gamma_r$. Combining
these local area bounds with the {\it a priori} curvature
estimates for minimizing surfaces, we can take a sequence of $r$'s
going to infinity and find a subsequence of $\Gamma_r$'s that
converge to a {\underline{complete}} area-minimizing surface
\begin{equation}
    \Gamma \subset \{ - h \leq x_3 \leq h \} \, .
\end{equation}
Since $\Gamma$ is pinched between the planes $\{ x_3 = \pm h\}$,
the estimates for minimizing surfaces implies that (outside a
large compact set) $\Gamma$ is a graph over the plane $\{ x_3 = 0
\}$ and hence has quadratic area growth and finite total
curvature.  Moreover, using the form of the metric $g_{ij}$, we
see that $|\nabla u|$ decays like $|x|^{-1}$ and
\begin{equation}    \label{e:kg1}
    \int_{\sigma_s} k_g  = (2 \, \pi \, s + O(1)) \, (s^{-1} + O(s^{-2}))
    = 2 \, \pi + O(s^{-1}) \, ,
\end{equation}
where $\sigma_s$ is the curve $ \{ x_1^2 + x_2^2 = s^2 \} \cap
\Gamma$ and $k_g$ is the geodesic curvature of $\sigma_s$
(considered as a curve in the surface $\Gamma$).

To get the contradiction, one combines stability of $\Gamma$ with
the positive scalar curvature of $M$ to see that no such $\Gamma$
could have existed.  Namely,  substituting the Gauss equation into
the stability inequality in a general $3$-manifold (see, e.g.,
\cite{CM1})
 gives
\begin{equation}    \label{e:stabipm}
\int_{\Gamma} (|A|^2 / 2+  \Scal_M -   K_{\Sigma} ) \phi^2 \leq \,
\int_{\Gamma} |\nabla \phi|^2  \, .
\end{equation}
Since  $\Gamma$ has quadratic area growth, we can choose a
sequence of (logarithmic) cutoff functions in \eqr{e:stabipm} to
get
\begin{equation}    \label{e:stabipm2}
0 < \int_{\Sigma} (|A|^2 / 2 +  \Scal_M )
      \leq  \int_{\Sigma}   K_{\Sigma} < \infty \, ;
\end{equation}
since $K_{\Sigma}$ may not be positive, we also used that $\Gamma$
has finite total curvature.  Moreover, we used that $\Scal_M$ is
positive outside a compact set to see that the first integral  in
\eqr{e:stabipm2} was positive.   Finally,  substituting
\eqr{e:stabipm2} into the Gauss-Bonnet formula gives that
$\int_{\sigma_s} k_g$ is {\underline{strictly}} less than $2 \pi$
for $s$ large, contradicting \eqr{e:kg1}.
\end{proof}

\subsection{Black holes}

Another way that minimal surfaces enter into relativity is through
black holes.   Suppose that we have a three-dimensional time-slice
$M$ in a $3+1$-dimensional space--time. For simplicity, assume
that $M$ is totally geodesic and hence has non-negative scalar
curvature.
  A  closed surface $\Sigma$ in $M$ is  said
to be trapped if its mean curvature is everywhere negative with
respect to its outward normal.  Physically, this means that the
surface emits an outward shell of light whose surface area is
decreasing everywhere on the surface.  The existence of a closed
trapped surface implies the existence of a black hole in the
space-time.

 Given a
trapped surface, we can look for the outermost trapped surface
containing it; this outermost surface is called an apparent
horizon.  It is not hard to see that an apparent horizon must be a
minimal surface and, moreover, a barrier argument shows that it
must be stable. Since $M$ has non-negative scalar curvature,
stability in turn implies that it must be diffeomorphic to a
sphere.  See, for instance, \cite{Br} and \cite{HuI} for some
results on black holes, horizons, etc.

\subsection{Constant mean curvature surfaces}

At least since the time of Plateau, minimal surfaces have been
used to model soap films.  This is because the mean curvature of
the surface models the surface tension and this is essentially the
only force acting on a soap film. Soap bubbles, on other hand,
enclose a volume and thus the pressure gives a second
counterbalancing force.  It follows easily that these two forces
are in equilibrium when the surface has constant mean curvature.

For the same reason,   constant mean curvature surfaces arise in
the isoperimetric problem.  Namely, a surface that minimizes
surface area while enclosing a fixed volume must have constant
mean curvature (or ``cmc'').  It is not hard to see that such an
isoperimetric surface in $\RR^n$ must be a round sphere.   There
are two interesting partial converses to this.  First, by a
theorem of H. Hopf, any cmc $2$-sphere in $\RR^3$   must be round.
Second, using the maximum principle (``the method of moving
planes'') A.D. Alexandrov showed that any closed embedded cmc
hypersurface in $\RR^n$ must be a round sphere. It turned out,
however, that not every closed immersed cmc surface is round.  The
first examples were immersed cmc tori constructed by H. Wente,
\cite{Wn}.   N. Kapouleas constructed many new examples, including
closed higher genus cmc surfaces, \cite{Ka1}, \cite{Ka2}.  See
also \cite{MzPaP} and \cite{KoKu} for other results on the space
of such cmc surfaces.

Many of the techniques developed for studying minimal surfaces
generalize to general constant mean curvature surfaces.

\subsection{Finite extinction for Ricci flow}

We close this survey with indicating how minimal surfaces can be
used to show that on a homotopy $3$-sphere the Ricci flow become
extinct in finite time (see \cite{CM12}, \cite{Pe} for details).

Let $M^3$ be a smooth closed orientable $3$--manifold and let
$g(t)$ be a one--parameter family of metrics on $M$ evolving by
the Ricci flow, so
\begin{equation}  \label{e:eqRic}
 \partial_t g=-2\,\Ric_{M_t}\, .
\end{equation}

In an earlier section, we saw that there is a natural way of
constructing minimal surfaces on many $3$-manifolds and that comes
from the min--max argument where the minimal of all maximal slices
of sweep--outs is a minimal surface. The idea is then to look at
how the area of this min--max surface changes under the flow.
Geometrically the area measures a kind of width of the
$3$--manifold and as we will see for certain $3$--manifolds
(those, like the $3$--sphere, whose prime decomposition contains
no aspherical factors) the area becomes zero in finite time
corresponding to that the solution becomes extinct in finite time.

\begin{figure}[htbp]
\begin{center}
    \input{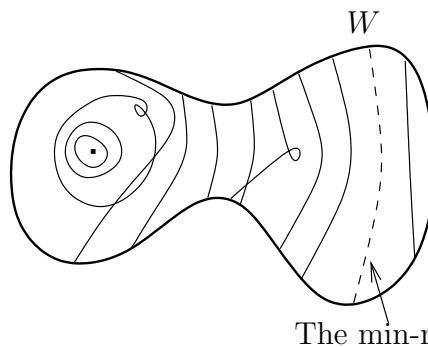}
    \caption{The sweep--out, the min--max surface, and the width W. First published in the Journal of the
    American Mathematical Society in 2005,
    published by the American Mathematical Society.}
    \label{f:1}
\end{center}
\end{figure}

Fix a continuous map $\beta: [0,1] \to C^0\cap L_1^2 (\SS^2 , M)$
where $\beta (0)$ and $\beta (1)$ are constant maps so that
$\beta$ is in  the nontrivial homotopy class $[\beta]$ (such
$\beta$ exists when $M$ is a homotopy $3$-sphere). We define the
width  $W=W(g,[\beta])$ by
\begin{equation}    \label{e:w3}
   W(g) = \min_{\gamma \in [ \beta]}
\, \max_{s \in [0,1]} \Energy (\gamma(s)) \,
   .
\end{equation}

\vskip2mm The next theorem gives an upper bound for the derivative
of $W(g(t))$ under the Ricci flow which forces the solution $g(t)$
to become extinct in finite.

\begin{theorem}     \label{t:upper}
Let $M^3$ be a homotopy $3$-sphere equipped with a Riemannian
metric $g=g(0)$. Under the Ricci flow, the width $W(g(t))$
satisfies
\begin{equation}   \label{e:di1a}
\frac{d}{dt} \, W(g(t))  \leq -4 \pi + \frac{3}{4 (t+C)} \,
W(g(t))   \, ,
\end{equation}
in the sense of the limsup of forward difference quotients. Hence,
 $g(t)$  must become extinct in finite time.
\end{theorem}

The $4\pi$ in \eqr{e:di1a} comes from the Gauss--Bonnet theorem
and the $3/4$ comes from the bound on the minimum of the scalar
curvature that the evolution equation implies.  Both of these
constants matter whereas the constant $C$ depends on the initial
metric and the actual value is not important.

To see that \eqr{e:di1a} implies finite extinction time rewrite
\eqr{e:di1a} as
\begin{equation}
\frac{d}{dt} \left( W(g(t)) \, (t+C)^{-3/4} \right) \leq - 4\pi \,
(t+C)^{-3/4}
\end{equation} and
integrate to get
\begin{equation}  \label{e:lastaa}
 (T+C)^{-3/4} \, W(g(T)) \leq C^{-3/4} \, W(g(0))
- 16 \, \pi \, \left[ (T+C)^{1/4} - C^{1/4} \right]   \, .
\end{equation}
Since $W \geq 0$ by definition and the right hand side of
\eqr{e:lastaa} would become negative for $T$ sufficiently large we
get the claim.

\vskip2mm As a corollary of this theorem we get finite extinction
time for the Ricci flow.

\begin{corollary}  \label{c:upper}
Let $M^3$ be a homotopy $3$-sphere equipped with a Riemannian
metric $g=g(0)$. Under the Ricci flow $g(t)$ must become extinct
in finite time.
\end{corollary}

\bibliographystyle{plain}
%\bibliography{bib__File}

\begin{thebibliography}{A}

\bibitem{Al}
W. Allard,   {\it On the first variation of a varifold},
 Ann. of Math. (2) 95 (1972) 417--491.
\bibitem{Am1}
F. J. Almgren,   Jr., {\it Some interior regularity theorems for
minimal surfaces and an extension of Bernstein's theorem}, Ann. of
Math. (2) 84 (1966) 277--292.
\bibitem{Am2}
F. J. Almgren,   Jr., {\it Almgren's big regularity paper.
$Q$-valued functions minimizing Dirichlet's integral and the
regularity of area-minimizing rectifiable currents up to
codimension 2}, With a preface by Jean E. Taylor and Vladimir
Scheffer. World Scientific Monograph Series in Mathematics, 1.
World Scientific Publishing Co., Inc., River Edge, NJ, 2000.


\bibitem{Am3}
F. J. Almgren,   Jr., {\it Minimal surface forms}, Math. Intell. 4
(1982) 164--172.

\bibitem{Be}
S. Bernstein,  {\it \"Uber ein geometrisches Theorem und seine
Anwendung auf die partiellen Differentialgleichungen vom
ellipschen} Typos. Math. Zeit. 26 (1927) 551--558 (translation of
the original version in Comm. Soc. Math. Kharkov 2-\`eme s\'er. 15
(1915--17) 38--45).
\bibitem{Ber}
L. Bers, {\it Isolated singularities of minimal surfaces},
 Ann. of Math. (2)  53 (1951) 364--386.

 \bibitem{BJO}
G.P. Bessa, L. Jorge and G. Oliveira-Filho,  Half-space theorems
for minimal surfaces with bounded curvature, {\it J. Diff. Geom.}
57 (2001) 493--508.

\bibitem{BDGG}
E. Bombieri, E. De Giorgi,  and E. Giusti,  {\it Minimal cones and
the Bernstein problem}, Invent. Math. 7 (1969) 243--268.

\bibitem{Br}
H. Bray, {\it Proof of the Riemannian Penrose inequality using the
positive mass theorem}, J. Differential Geom. 59 (2001), no. 2,
177--267.


\bibitem{Ca}
E. Calabi, Problems in differential geometry, Ed. S. Kobayashi and
J. Eells, Jr., {\it Proceedings of the United States-Japan Seminar
in Differential Geometry, Kyoto, Japan, 1965}. Nippon Hyoronsha
Co., Ltd., Tokyo (1966) 170.
\bibitem{Ch}
S.S. Chern,   {\it The geometry of $G$-structures},  Bull. Amer.
Math. Soc. 72 (1966) 167--219.
\bibitem{CiSc}
H.I. Choi and R. Schoen, {\it The space of minimal embeddings of a
surface into a three-dimensional manifold of positive Ricci
curvature},   Invent. Math.  81 (1985) 387--394.


\bibitem{CD}
T.H. Colding and C. De Lellis, {\it The min--max construction of
minimal surfaces}, Surveys in differential geometry, Vol. 8,
Lectures on Geometry and Topology held in honor of Calabi, Lawson,
Siu, and Uhlenbeck at Harvard University, May 3--5, 2002,
Sponsored by the Journal of Differential Geometry, (2003) 75--107,
math.AP/0303305.



\bibitem{CM1}
T.H. Colding and W.P. Minicozzi II, \textit{Minimal surfaces,
Courant Lecture Notes in Math.}, v. 4, 1999.
\bibitem{CM2}
T.H. Colding and W.P. Minicozzi II, \textit{Estimates for
parametric elliptic integrands},    International Mathematics
Research Notices, no. 6 (2002) 291-297.
\bibitem{CM3}
T.H. Colding and W.P. Minicozzi II, {\it The space of embedded
minimal surfaces of fixed genus in a $3$-manifold I; Estimates off
the axis for disks},  Annals of Math., 160 (2004) 27--68,
math.AP/0210106.
\bibitem{CM4}
T.H. Colding and W.P. Minicozzi II, {\it The space of embedded
minimal surfaces of fixed genus in a $3$-manifold II; Multi-valued
graphs in disks},  Annals of Math., 160 (2004) 69--92,
math.AP/0210086.
\bibitem{CM5}
T.H. Colding and W.P. Minicozzi II, {\it The space of embedded
minimal surfaces of fixed genus in a $3$-manifold III; Planar
domains},  Annals of Math., 160 (2004) 523--572, math.AP/0210141.
\bibitem{CM6}
T.H. Colding and W.P. Minicozzi II, {\it The space of embedded
minimal surfaces of fixed genus in a $3$-manifold IV; Locally
simply connected},  Annals of Math., 160 (2004) 573--615,
math.AP/0210119.
\bibitem{CM7}
T.H. Colding and W.P. Minicozzi II, {\it The space of embedded
minimal surfaces of fixed genus in a $3$-manifold V; Fixed genus},
math.DG/0509647.
\bibitem{CM8}
T.H. Colding and W.P. Minicozzi II, \textit{Multi-valued minimal
graphs and properness of disks},    International Mathematics
Research Notices, no. 21 (2002) 1111-1127.
\bibitem{CM9}
T.H. Colding and W.P. Minicozzi II, \textit{On the structure of
embedded minimal annuli},   International Mathematics Research
Notices, no. 29 (2002) 1539--1552.


\bibitem{CM10}
T.H. Colding and W.P. Minicozzi II, \textit{Disks that are double
spiral staircases}, Notices of the AMS, Vol. 50, no. 3, March
(2003) 327--339.

\bibitem{CM11}
T.H. Colding and W.P. Minicozzi II, \textit{Embedded minimal
disks: Proper versus nonproper - global versus local},
{Transactions of the AMS}, 356 (2004) 283-289, math.DG/0210328.
\bibitem{CM12}
T.H. Colding and W.P. Minicozzi II, \textit{Estimates for the
extinction time for the Ricci flow on certain $3$--manifolds and a
question of Perelman}, JAMS,  18  (2005),  no. 3, 561--569,
math.AP/0308090.
\bibitem{CM13}
T.H. Colding and W.P. Minicozzi II, \textit{ Examples of embedded
minimal tori without area bounds},    International Mathematics
Research Notices, 99 no. 20 (1999) 1097--1100.
\bibitem{CM14}
T.H. Colding and W.P. Minicozzi II,  \textit{ Complete properly
embedded minimal surfaces in $\RR^3$}, Duke Math. J. 107 (2001)
421--426.
\bibitem{CM15}
T.H. Colding and W.P. Minicozzi II,  \textit{ Embedded minimal
disks}, Global theory of minimal surfaces, 405--438, Clay Math.
Proc., 2, Amer. Math. Soc., Providence, RI, 2005, math.DG/0206146.

\bibitem{CM16}
T.H. Colding and W.P. Minicozzi II, {\it The Calabi--Yau
conjectures for embedded surfaces}, preprint, math.DG/0404197.

\bibitem{Co}
P. Collin,   {\it Topologie et courbure des surfaces minimales
proprement plong‰es de $\RR^3$},   Ann. of Math. (2) 145 (1997)
1--31.

\bibitem{DG}
E. De Giorgi, {\it Frontiere orientate di misura minima}, Sem.
Mat. Scuola Norm. Sup. Pisa (1961) 1--56.
\bibitem{De}
B. Dean, {\it Compact Embedded Minimal Surfaces of Positive Genus
Without Area Bounds}, Geom. Ded. 102 (2003), 45--52.

\bibitem{EWWi}
T. Ekholm, B. White, and D. Wienholtz, {\it Embeddedness of
minimal surfaces with total boundary curvature at most $4\pi$},
 Ann. of Math. (2) 155 (2002), no. 1,
209--234.


\bibitem{Fe}
H. Federer, {\it  Geometric measure theory}, Springer-Verlag,
Berlin--Heidelberg--New York, 1969.

\bibitem{Fi}
R. Finn, {\it Capillary surface interfaces}, Notices Amer. Math.
Soc. 46 (1999), no. 7, 770--781.

\bibitem{FHS}
M.H. Freedman, J. Hass, and P. Scott, {\it Least area
incompressible surfaces in $3$-manifolds}, Invent. Math. 71
(1983), no. 3, 609--642.
\bibitem{Fr}
A. Fraser, {\it On the free boundary variational problem for
minimal disks}, Comm. Pure Appl. Math. 53 (2000) 931--971.


\bibitem{GuSp}
R. Gulliver and J. Spruck,  {\it On embedded minimal surfaces},
Ann. of Math. (2) 103 (1976) 331--347; Ann. of Math. (2) 109
(1979) 407--412.

\bibitem{HaKLi}
R. Hardt, D. Kinderlehrer, and F. H. Lin, {\it The variety of
configurations of static liquid crystals}, Variational methods
(Paris, 1988), 115--131, Progr. Nonlinear Differential Equations
Appl., 4, Birkhuser Boston, Boston, MA, 1990.


\bibitem{HaSi}
R. Hardt and L. Simon, {\it Boundary regularity and embedded
solutions for the oriented Plateau problem}, Ann. of Math. (2) 110
(1979), no. 3, 439--486.

\bibitem{He}
E. Heinz,  {\it \"Uber die L\"osungen der
Minimalfl\"achengleichung},
 Nachr. Akad. Wiss. G\"ottingen Math.--Phys. Kl, II (1952) 51--56.

\bibitem{Ho}
D. Hoffman, \textit{Mixing materials and mathematics},  Nature
384 (1996) 28 - 29.


\bibitem{HoK}
D. Hoffman and H. Karcher, \textit{Complete embedded minimal
surfaces with finite total curvature}, Geometry V (R. Osserman,
ed.) Encyclopaedia Math. Sci. 90, Springer-Verlag, New York (1997)
5--93.
\bibitem{HoMe}
D. Hoffman and W. Meeks  III,   {\it A complete embedded minimal
surface in $\RR^3$ with genus one and three ends},
 J. Diff. Geom. 21 (1985) 109--127.

\bibitem{HoWeWo1}
D. Hoffman, M. Weber, and M. Wolf,   {\it   An embedded genus-one
helicoid}, math.DG/0401080.

\bibitem{HoWeWo2}
D. Hoffman, M. Weber, and M. Wolf,   {\it   An embedded genus-one
helicoid}, PNAS, November 15, 2005, vol. 102, no. 46.

\bibitem{HuI} G. Huisken and T. Ilmanen, {\it The inverse mean
curvature flow and the Riemannian Penrose inequality}, J.
Differential Geom. 59 (2001), no. 3, 353--437.

\bibitem{I}
T. Ilmanen, {\it A strong maximum principle for singular minimal
hypersurfaces},
 Calc. Var. Partial Differential Equations 4 (1996), no. 5, 443--467.

\bibitem{JXa1}
L. Jorge and F. Xavier, {\it On the existence of complete bounded
minimal surfaces in $ \RR^{n}$},  Bol. Soc. Brasil. Mat. 10
(1979), no. 2, 171--173.
\bibitem{JXa2}
L. Jorge and F. Xavier, {\it A complete minimal surface in $\RR^3$
between two parallel planes},  Annals of Math. (2) 112 (1980)
203--206.
\bibitem{JXa3}
L. Jorge and F. Xavier, {\it An inequality between the exterior
diameter and the mean curvature of bounded immersions},   Math.
Zeit. 178 (1981), no. 1, 77--82.


\bibitem{Ka1}
N. Kapouleas, \textit{Compact constant mean curvature surfaces in
Euclidean three-space}, J. Differential Geom. 33 (1991), no. 3,
683--715.


\bibitem{Ka2}
N. Kapouleas, \textit{Complete constant mean curvature surfaces in
Euclidean three-space}, Ann. of Math. (2) 131 (1990), no. 2,
239--330.


\bibitem{Ka3}
N. Kapouleas, \textit{Complete embedded minimal surfaces of finite
total curvature}, J. Diff. Geom., 47 (1997) 95--169.


\bibitem{KoKu}
N. Korevaar and R. Kusner, {\it The global structure of constant
mean curvature surfaces}, Invent. Math. 114 (1993), no. 2,
311--332.


\bibitem{La}
H. B. Lawson,  {\it Lectures on minimal submanifolds}, vol. I,
Publish or Perish, Inc, Berkeley, 1980.


\bibitem{LMaMo1}
F. Lopez, F.  Martin, and S. Morales,  {\it Adding handles to
Nadirashvili's surfaces}, J. Diff. Geom. 60 (2002), no. 1,
155--175.
\bibitem{LMaMo2}
F. Lopez, F.  Martin, and S. Morales, {\it Complete nonorientable
minimal surfaces in a ball of $\RR^3$}, preprint.

\bibitem{LoRo}
F. Lopez and A. Ros, {\textit{On embedded complete minimal
surfaces of genus zero}}, J. Differential Geom. 33 (1991), no. 1,
293--300.

\bibitem{MaMo1}
F.  Martin and S. Morales,   {\it A complete bounded minimal
cylinder in $\RR^3$},  Michigan Math. J. 47 (2000), no. 3,
499--514.
\bibitem{MaMo2}
F.  Martin and S. Morales,  {\it On the asymptotic behavior of a
complete bounded minimal surface in $\RR^3$},    Trans. Amer.
Math. Soc. 356 (2004), 3985--3994.
\bibitem{MaMo3}
F.  Martin and S. Morales, {\it Complete proper minimal surfaces
in convex bodies of $\RR^3$}, Duke Math. J. 128 (2005), no. 3,
559--593.


\bibitem{MzPaP}
R. Mazzeo, F. Pacard, and D. Pollack,{\it Connected sums of
constant mean curvature surfaces in Euclidean 3 space}, J. Reine
Angew. Math. 536 (2001), 115--165.

\bibitem{Me1}
W. Meeks III,  {\it The regularity of the singular set in the
Colding and Minicozzi lamination theorem}, Duke Math. Jour., 123
(2004), no. 2, 329--334.
\bibitem{Me2}
W. Meeks III,    {\it The lamination metric for the
Colding-Minicozzi minimal lamination},   Illinois J. Math. 49
(2005), no. 2, 645--658.


\bibitem{MeP}
W. Meeks III and J. Perez,   {\it Conformal properties in
classical minimal surface theory}, Surveys in Diff. Geom. IX:
Eigenvalues of Laplacians and other geometric operators, Ed. by A.
Grigor'yan and S.T. Yau, International Press (2004), 275--335.
\bibitem{MePRs1}
W. Meeks III, J. Perez, and A. Ros, {\it The geometry of minimal
surfaces of finite genus I; Curvature estimates and
quasiperiodicity}, J. Differential Geom. 66 (2004), 1--45.
\bibitem{MePRs2}
W. Meeks III, J. Perez, and A. Ros, {\it The geometry of minimal
surfaces of finite genus II; Nonexistence of one limit end
examples}, Invent. Math. 158 (2004), no. 2, 323--341.
\bibitem{MePRs3}
W. Meeks III, J. Perez, and A. Ros, {\it The geometry of minimal
surfaces of finite genus III; bounds on the topology and index of
classical minimal surfaces}, preprint.
\bibitem{MeRo1}
W. Meeks III and H. Rosenberg, {\it The uniqueness of the helicoid
and the asymptotic geometry of properly embedded minimal surfaces
with finite topology}, Ann. of Math.,  161  (2005),  no. 2,
727--758.

\bibitem{MeRo2}
W. Meeks III and H. Rosenberg, {\it The  minimal lamination
closure theorem}, preprint.

\bibitem{MeSiYa}
W. Meeks  III, L. Simon, and S.T. Yau,   {\it Embedded minimal
surfaces, exotic spheres and manifolds with positive Ricci
curvature},
 Ann. of Math. (2) 116 (1982) 621--659.
 \bibitem{MeWe}
W. Meeks III and M. Weber, {\it Existence of bent helicoids and
the regularity of the singular set in the Colding-Minicozzi
lamination theorem}, in preparation.
\bibitem{MeYa1}
W. Meeks  III and S.T. Yau,   {\it The classical Plateau problem
and the topology of three dimensional manifolds},
 Topology 21 (1982) 409--442.
\bibitem{MeYa2}
W.H. Meeks and S.T. Yau,  \textit{Topology of three--dimensional
manifolds and the embedding problems in minimal surface theory},
Ann. of Math. (2) 112 (1980), no. 3, 441--484.
 \bibitem{MiMo}
M.J. Micallef and J.D. Moore, {\it Minimal two--spheres and the
topology of manifolds with positive curvature on totally isotropic
two--planes}, Ann. of Math. (2) 127 (1988) no. 1 199--227.

\bibitem{Na1}
N. Nadirashvili, {\it Hadamard's and Calabi-Yau's conjectures on
negatively curved and minimal surfaces},  Invent. Math. 126 (1996)
457--465.
\bibitem{Na2}
N. Nadirashvili,  {\it An application of potential analysis to
minimal surfaces},  Mosc. Math. J. 1 (2001), no. 4, 601--604, 645.

\bibitem{Os}
R. Osserman, {\it A survey of minimal surfaces},   Dover, 2nd.
edition (1986).


\bibitem{Pe}
G. Perelman, {\it Finite extinction time for the solutions to the
Ricci flow on certain three--manifolds},  math.DG/0307245.

\bibitem{Pz}
J. Perez, {\it Limits by rescalings of minimal surfaces: Minimal
laminations, curvature decay and local pictures, notes for  the
workshop "Moduli Spaces of Properly Embedded Minimal Surfaces"},
American Institute of Mathematics, Palo Alto, California (2005).

\bibitem{Pi}
J.T. Pitts, {\it Existence and regularity of minimal surfaces on
Riemannian manifolds}, Princeton University Press, Princeton,
N.J.; University of Tokyo Press, Tokyo (1981).
\bibitem{PiRu}
J.T. Pitts and J.H. Rubinstein, {\it Applications of minmax to
minimal surfaces and the topology of $3$--manifolds},
Miniconference on geometry and partial differential equations,
Proceedings of the CMA, Australia National University   (1986).


\bibitem{Ro1}
H. Rosenberg, {\it Some recent developments in the theory of
properly embedded minimal surfaces in $\RR^3$}, Seminare Bourbaki
1991/92, Asterisque No. 206 (1992) 463--535.

\bibitem{Ro2}
H. Rosenberg,   {\it Some recent developments in the theory of
minimal surfaces in 3-manifolds},  IMPA Mathematical Publications.
24th Brazilian Mathematics Colloquium, Instituto de Matemática
Pura e Aplicada (IMPA), Rio de Janeiro, 2003.

\bibitem{Ro3}
H. Rosenberg, A complete embedded minimal surface in $\RR^3$ of
bounded curvature is proper, preprint.


\bibitem{SaUh}
J. Sacks and K. Uhlenbeck, \textit{The existence of minimal
immersions of $2$--spheres}, Ann. of Math. (2) 113 (1981) no. 1,
1--24.

\bibitem{Sc1}
R. Schoen,  {\it Estimates for stable minimal surfaces in
three--dimensional manifolds},
 In Seminar on Minimal Submanifolds,
Ann. of Math. Studies, vol. 103, Princeton University Press,
Princeton, N.J., (1983)  111--126.
\bibitem{Sc2}
 R.  Schoen, {\it Uniqueness, symmetry, and embeddedness of minimal surfaces},
  J. Differential Geom. 18 (1983), no. 4, 791--809 (1984).

\bibitem{ScSi}
R. Schoen and L.  Simon,   {\it Regularity of simply connected
surfaces with quasi-conformal Gauss map}, In Seminar on Minimal
Submanifolds,
 Annals of Math. Studies, vol. 103,  Princeton University Press,
Princeton, N.J., (1983) 127--145.
\bibitem{ScSiYa}
R. Schoen, L. Simon,   and S.T. Yau,   {\it Curvature estimates
for minimal hypersurfaces},
 Acta Math. 134 (1975) 275--288.

 \bibitem{ScYa1}
Schoen, R.  and Yau, S. T., {\it On the proof of the positive mass
conjecture in general relativity}, Comm. Math. Phys. 65 (1979),
no. 1, 45--76.
 \bibitem{ScYa2}
R. Schoen and S.T. Yau, {\it Existence of incompressible minimal
surfaces and the topology of three dimensional manifolds with
nonnegative scalar curvature}, Ann. of Math. (2) 110 (1979)
127--142.



\bibitem{Si1}
L. Simon, {\it Asymptotic behaviour of minimal graphs over
exterior domains},  Ann. Inst. H. Poincaré Anal. Non Linéaire 4
(1987)  231--242.
\bibitem{Si2}
  L. Simon,   {\it Remarks on curvature
estimates for minimal hypersurfaces},
 Duke Math. J. 43 (1976) 545--553.
\bibitem{Si3}
  L. Simon,   {\it Singularities of Geometric Variational Problems}, In Nonlinear Partial
Differential Equations in Differential Geometry (R. Hardt and M.
Wolf, Ed.),
   American Mathematical Society, Providence  (1996) 185--223.

\bibitem{Si4}
  L. Simon,
 {\it A strict maximum principle for area minimizing
hypersurfaces}, J. Differential Geom. 26 (1987), no. 2, 327--335.



\bibitem{Sim}
J. Simons,   {\it Minimal varieties in Riemannian manifolds},
 Ann. of Math. (2) 88
 (1968) 62--105.
 \bibitem{Sm}
F. Smith, {\it On the existence of embedded minimal $2$--spheres
in the $3$--sphere, endowed with an arbitrary Riemannian metric},
supervisor L. Simon, University of Melbourne (1982).

\bibitem{SoWh} B. Solomon and B. White,  {\it A strong maximum
principle for varifolds that are stationary with respect to even
parametric elliptic functionals}, Indiana Univ. Math. J. 38
(1989), no. 3, 683--691.

\bibitem{Ta}
J. Taylor, {\it Some mathematical challenges in materials
science}, Mathematical challenges of the 21st century (Los
Angeles, CA, 2000). Bull. Amer. Math. Soc. (N.S.) 40 (2003), no.
1, 69--87.

\bibitem{Th}
D.W. Thompson, {\it On growth and form}, New ed., Cambridge Univ.
Press, Cambridge (1942).

\bibitem{Tr}
M. Traizet,   {\it Adding handles to Riemann's minimal surfaces},
 J. Inst. Math. Jussieu 1
 (2002) 145--174.
\bibitem{WeWo}
M. Weber and M. Wolf, {\it  Teichm\"uller theory and handle
addition for minimal surfaces}, Ann. of Math. (2), 156 (2002)
713--795.

\bibitem{Wn}
H. Wente,   {\it Counterexample to a conjecture of H. Hopf},
Pacific J. Math. 121 (1986), no. 1, 193--243.


 \bibitem{W}
B. White,   {\it Evolution of curves and surfaces by mean
curvature}, Proceedings of the ICM  (2002).

\bibitem{Xa}
F. Xavier, {\it Convex hulls of complete minimal surfaces},
   Math. Ann. 269 (1984) 179--182.

 \bibitem{Y}
S.D. Yang, {\it A connected sum construction for complete minimal
surfaces of finite total curvature}, Comm. Anal. Geom. 9 (2001),
no. 1, 115--167.
 \bibitem{Ya1} S.T. Yau, {\it Nonlinear analysis in geometry},
L'Eseignement Mathematique (2) 33 (1987) 109--158.


\end{thebibliography}

\end{document}